%% The Frobenius structure of local cohomology  
%% AMSTEX, amsppt documentstyle

\input amstex
     %%%   Version submitted for next evaluation by referee
\documentstyle{amsppt}
\magnification=1200
\hoffset .6truein
\voffset .8truein
\loadeurm
\loadeusm
\loadbold
\font\twbf=cmbx10 scaled\magstep2
\font\chbf=cmbx10 scaled\magstep1
\font\hdbf=cmbx10
 
                     %%% bold math italic
\font\smrm=cmr7     
\pagewidth{6.5 true in}
\pageheight{8.9 true in}
\nologo
\overfullrule=3pt
\define\0{^{@,@,\circ}}  
\define\1{^{-1}}
\define\2{^{\hbox{\rm h}}}
\define\8{^{\infty}} 
\define\({\bigl(}
\define\){\bigr)}
\define\<{\langle}
\define\>{\rangle}

\define\an{anti\hskip .15pt -nilpotent}
\define\Ann{\hbox{\rm Ann}}

\define\blank{\underline{\phantom{J}}}
\define\bpf{\demb{Proof}}

\define\cH{{\Cal H}}

\define\clo#1#2#3{{{#1}^{#2}}_{#3}}

\define\col#1{\,\colon_{#1}}

\define\dl{{\Delta}}
\define\demb#1{\demo{\bf #1}}

\define\di{\text{dim}\,}
\define\dlim#1{\underset{\longrightarrow}\to{\hbox{\rm lim}}\,_{#1}\,}

\define\Ext{\hbox{\rm Ext}}
\define\f{\{\hskip -.55pt F \hskip -.4pt\}}
\define\F{\boldkey F}
\define\fA{{\frak A}}
\define\fh{finite FH-length}

\define\Hom{\hbox{\rm Hom}}

\define\id{\hbox{\rm id}}
\define\inx{\langle x^{-1}\rangle}

\define\im{\hbox{\rm Im}\,}
\define\imp{\Rightarrow}
\define\inc{\subseteq}
\define\inj{\hookrightarrow}
\define\intp #1 {\lfloor #1 \rfloor}
\define\link{\hbox{\rm link}}

\define\Ker{\text{Ker}\,}

\define\n{\frak{n}}
\define\N{{\Bbb N}}
\define\ng{\hbox{\rm neg}}
\define\ngs{\hbox{\smrm neg}}
\define\noi{\noindent}

\define\p{^{[p]}}
\def\pa{\parindent 35 pt}
\define\part#1{\item{\hbox{\rm (#1)}}}

\define\Pnl{\boldkey{(}}
\define\Pnr{\boldkey{)}}
\define\pos{\hbox{\rm pos}}

\define\ps{{}^\Psi\!}
\define\ov#1{\overline{#1}}
\define\q{^{[q]}}
\define\qedd{\qed\enddemo}

\define\R{\widehat{R}}
\define\Rad{\text{Rad}\,}

\define\rf{R\{\hskip -.55pt F \hskip -.4pt\}}
\define\rmk{(R,\,m,\,K)}
\def\sy#1 {^{\Pnl#1\Pnr}}

\define\spu{\hbox{\rm supp}}
\define\spus{\hbox{\smrm supp}}
\define\surj{\twoheadrightarrow}
\define\surjb{\twoheadleftarrow}

\define\ub{^{\bullet}}
\define\uf{\underline{f}}
\define\usurj{\lower .3 pt \hbox{$\uparrow$} \kern - 5.1pt \hbox{\raise 1.1 pt \hbox{$\uparrow$}} }
\define\ux{\underline{x}}

\define\vect#1#2{#1_1,\, \ldots, #1_{#2}}
\define\w{\omega}
\define\wh{\widehat}
\define\wt{\widetilde}
\define\Z{{\Bbb Z}}
\quad\bigskip\bigskip\bigskip
\centerline{\twbf THE FROBENIUS STRUCTURE} 
\vskip 10 pt plus .5 pt minus .5 pt
\centerline{\twbf OF LOCAL COHOMOLOGY}

\topmatter
\author
by Florian Enescu and Melvin Hochster
\endauthor
\rightheadtext{The Frobenius Structure of Local Cohomology}
\leftheadtext{Florian Enescu and Melvin Hochster}

\abstract Given a local ring of positive prime characteristic there is a natural Frobenius action on its 
local cohomology modules with support at its maximal ideal.
In this paper we study the local rings for which the local cohomology modules have only finitely 
many submodules invariant under the Frobenius action. In 
particular we prove that F-pure Gorenstein local rings as well as the face 
ring of a finite simplicial complex localized or completed at its 
homogeneous maximal ideal have this property. We also introduce the notion of an 
anti-nilpotent Frobenius action on an Artinian module
over a local ring and use it to study those rings for which the lattice of submodules of 
the local cohomology that are invariant under
Frobenius satisfies the Ascending Chain Condition. 
\endabstract

\thanks{The first author was partially supported by NSA Young Investigator grant H98230-07-1-0034;
the second author was partially supported by NSF grant DMS--0400633}
\endthanks
\endthanks
\thanks{Subject classification (MSC2000): 13A35, 13D45.}\endthanks
\endtopmatter

\document
{\baselineskip = 16 pt     %% Put this in once done
\parindent = 12 pt
\quad\bigskip\bigskip
{\chbf 
\centerline{\hdbf 1. INTRODUCTION} 
}
\bigskip

All given rings in this paper are commutative, associative with
identity, and Noetherian.  Throughout, $p$ denotes a positive
prime integer.  For the most part, we shall be studying local
rings, i.e., Noetherian rings with a unique maximal ideal. Likewise
our main interest is in rings of positive prime 
characteristic $p$.  If $(R,\,m) $ is local of characteristic $p$, there is
a natural action of the Frobenius endomorphism of $R$ on
each of its local cohomology modules $H^j_m(R)$.  
We call an $R$-submodule $N$ of one of these local cohomology
modules F-{\it stable} if the action of $F$ maps $N$ into itself. 

One of our objectives is to understand
when a local ring, $(R,\,m)$,  especially a reduced Cohen-Macaulay local ring,
has the property that only finitely many $R$-submodules of its local
cohomology modules F-stable.  When this occurs we say that $R$ is
FH-{\it finite}.    We shall also study the problem of determining conditions under
which the  local cohomology modules of $R$ have finite length in the category 
of $R$-modules with Frobenius action.  We say that $R$ has {\it finite FH-length}
in this case.
Of course, when the ring is Cohen-Macaulay there
is only one non-vanishing local cohomology module,  $H^d_m(R)$,
where $d = \di(R)$.  The problem of studying 
the F-stable submodules of $H^d_m(R)$ arises naturally
in tight closure theory, taking a point of view pioneered by K.\ Smith [Sm1,2,3,4]. 
E.g., if $R$ is complete, reduced, and Gorenstein,
the largest proper F-stable submodule of $H^d_m(R)$ corresponds to
the tight closure of 0 (in the finitistic sense:  see [HH2], \S8),
and its annihilator is the test ideal of $R$.  Also see
Discussion 2.10 here.  We would like to note here that other results 
related to F-stable submodules of local cohomology may be found in [En1], [En2], [Ka], [Sh].

The main result of Section 3 is one of general interest. Let $M$ be a
module over an excellent local ring $R$ and consider a family of submodules $\{N_\lambda\}_{\lambda \in \Lambda}$ 
of $M$ closed under sum, intersection and primary decomposition. Our result states, in
particular, that  if the set 
$\{ \Ann_R(M/N_\lambda): \lambda \in \Lambda\}$ consists of radical
ideals then it is finite. 

This theorem has a number of important corollaries. One of them, relevant to our objectives and also proven by Sharp 
(Theorem 3.10 and Corollary 3.11 in [Sh]), states that for any local ring $R$ of prime characteristic
$p>0$, if Frobenius acts injectively on an Artinian $R$-module $M$, then
the set of annihilators of F-stable submodules of $M$ is a finite
set of radical ideals closed under primary decomposition. This leads to the fact that if $R$ is F-pure and Gorenstein
(or even quasi-Gorenstein) then  $H^d_m(R)$ has only finitely F-stable  $R$-submodules. See \S 3. Another one is the 
fact that in an excellent local ring a family of radical ideals closed under sum, intersection and primary decomposition
is finite.  

In \S4 we explain the relationship between FH-finite rings and rings that have \fh. We introduce the notion of an {\it anti-nilpotent} Frobenius
action on an Artinian module over a local ring.  Using results of [Ly] and
[Ho4],  we show that the local cohomology of a local ring $R$ of characteristic
$p$ is anti-nilpotent if and only if the local cohomology of  $R[[x]]$ has
\fh, in which case the local cohomology of $R$ and every formal power series
ring over $R$ is FH-finite.

In \S5 we show that if $R$ is the face ring of a finite simplical complex localized
or completed at its homogeneous maximal ideal, then $R$ is FH-finite. See 
Theorem (5.1).

\quad\bigskip\bigskip
{\chbf 
\centerline{\hdbf 2. NOTATION AND TERMINOLOGY} 
}
\bigskip

\demb{Discussion 2.1:  Some basics about tight closure}
Unless otherwise specified,
we shall assume throughout that $R$ is a Noetherian ring of positive
prime characteristic $p$, although this hypothesis is usually
repeated in theorems and definitions.  $R\0$ denotes the complement of
the union of the minimal primes of $R$, and so, if $R$ is reduced,
$R\0$ is simply the multiplicative system of all nonzerodivisors in $R$.
We shall write $\F^e$ (or $\F^e_R$ if we need to specify the base ring)
for the {\it Peskine-Szpiro} or {\it Frobenius} functor from $R$-modules
to $R$-modules.  
Note that $\F^e$ preserves both freeness and
finite generation of modules, and is exact precisely when
$R$ is regular (cf.\ [Her], [Ku1]). If $N \inc M$  we write $N\q$ for 
the image of $F^e(N)$ in $F^e(M)$, although it depends on the inclusion
$N \to M$,  not just on $N$.  If $u \in M$ we write $u^{p^e}$
for the image $1 \otimes u$ of $u$ in $F^e(M)$. With this
notation,  $(u+v)^q = u^q + v^q$ and $(ru)^q = r^qu^q$ for
$u,\,v \in M$ and $r \in R$. 

From time to time, we assume some familiarity with basic tight closure theory
in prime characteristic $p > 0$.  We use the standard notation $N^*_M$ for
the tight closure of the submodule $N$ in the module $M$.  If $M$ is understood,
the subscript is omitted, which is frequently the case when $M = R$ and $N = I$
is an ideal.  We refer the reader to [HH1,2,3,6] and [Hu] for background in this area.

In particular, we assume cognizance of certain facts about test elements,
including the notion of a {\it completely stable} test element.  
We refer the reader to [HH2, \S6 and \S8 ],
[HH1], [HH3, \S6], and [AHH, \S2] for more information
about test elements and to \S3 of [AHH] for a discussion of
several basic issues related to the localization problem for
tight closure.    

\demb{Discussion 2.2: Local cohomology and the action of the 
Frobenius endomorphism} Our basic reference for local cohomology
is [GrHa]. Let $R$ be an arbitrary Noetherian ring, let $I$ be an
ideal of $R$ and let $M$ be any $R$-module.  The $i\,$th {\it local
cohomology module} $H^i_I(M)$ {with support in} $I$ may be obtained
in several ways.  It may be defined as $\dlim t \Ext^i_R(R/I^t,\,M)$:
here, any sequence of ideals cofinal with the powers of $I$ may
be used instead of the the sequence of powers, $\{I^t\}_t$.  
Alternatively, we may define $C^\bullet(f;\,R)$ to be the
complex $0 \to C^0 \to C^1 \to 0$  where  $C^0 = R$, $C^1 = R_f$
and the map is the canonical map $R \to R_f$,  and then if
$\uf$ is a sequence of elements $\vect f n$ we may define
$C^\bullet(\uf;\, R)$ to be the tensor product over $R$ of the $n$ complexes
$C^\bullet(f_i;\, R)$.  Finally, let $C^\bullet(\uf; M)$ denote
$C^\bullet(\uf; \,R) \otimes_R M$,  which has the form:
$$
0 \to M \to \bigoplus_i M_{f_i} \to \bigoplus_{i<j} M_{f_if_j} \to \cdots
\to M_{f_1\cdots f_n} \to 0 .
$$
The cohomology of this complex turns out to be $H^\bullet_I(M)$,
where $I = (\vect f n)R$, and actually depends only on
the radical of the ideal $I$.

By the standard theory of local duality (cf.\ [GrHa, Theorem (6.3)]) 
when  $(S,m_S,L)$ is Gorenstein with $\di(S) = n$ 
and $M$ is a finitely generated $S$-module,
$H^i_m(M) \cong \Ext^{n-i}_S(S,\,M)^\lor$ as functors of $M$,
where $N^\lor = \Hom_R\bigl(N,\,E_S(L)\bigr)$.  Here, $E_S(L)$
is an injective hull of $L$ over $S$. In particular,
if $\rmk$ is local of Krull dimension $d$ and is a homomorphic 
image of a Gorenstein local ring $S$ of dimension $n$,  then 
$\omega_R = \Ext^{n-d}(R,\,S)$ whose Matlis dual over $S$,
and, hence, over $R$ as well, is $H^d_m(R)$.  We refer to a finitely
generated $R$-module $\omega_R$
as a {\it canonical module} for $R$ if $\omega_R^\lor = H^d_m(R)$.  It is 
unique up to isomorphism, since its completion is dual to $H^d_m(R)$.
Our discussion shows that a canonical module exists if $R$ is a
homomorphic image of a Gorenstein ring;  in particular, $\omega_R$
exists if $R$ is complete.   When $R$ is Cohen-Macaulay, one 
has that $H^i_m(M) \cong \Ext^{d-i}_R(M,\, \omega_R)$ functorially for all 
finitely generated $R$-modules $M$.  

When $R$ is a normal local domain, $\omega_R$ is isomorphic
with as an $R$-module with an ideal of pure height one, i.e., with
a divisorial ideal.  

Finally, suppose that  $\rmk \to (S,\, m_S\, L)$ is a local homomorphism
such that $S$ a module-finite extension of $R$.  Let $\omega = 
\omega_R$ be a canonical module for $R$.  Then 
$\Hom_R(S,\, \omega)$ is a canonical module for $S$. Here, the
rings need not be Cohen-Macaulay, nor domains.  To see
this, note one can reduce at once to the complete case.  We have  
$$H^d_{m_S}(S) \cong  H^d_m(S) \cong S \otimes H^d_m(R).$$  Then
$E_S(L)$ may be identified with $\Hom_R\bigl(S,\, E_R(K)\bigr)$; 
moreover, on $S$-modules, the functors $\Hom_R\bigl(\blank, \, E_R(K)\bigr)$ and 
$\Hom_S(\blank,\, E_S(L))$ are isomorphic.  Hence
$$\Hom_S\big(H^d_{m_S}(S),\, E_S(L)\bigr) 
\cong \Hom_R(S \otimes_R H^d_m(R),\, E_R(K)\bigr).$$  By the
adjointness of tensor and Hom, this becomes
$$\Hom_R\bigl(S, \, \Hom_R\bigl(H^d_m(R), \, E_R(K)\bigr)\bigr) \cong
\Hom_R(S, \, \omega),$$  as required.

When $M = R$ we have an action of the Frobenius endomorphism
on the complex $C^\bullet(\uf;\,R)$ induced by the Frobenius
endomorphisms of the various rings $R_g$ where $g$ is a product
of a subset of $\vect f n$,  and the action on the cohomology
is independent of the choice of $f_i$.  

An alternative point of view is that, quite generally, if
$M \to M'$ is any map of $R$-modules then there is an
induced map  $H^i_I(M) \to H^i_I(M')$.  When  $S$ is an $R$-algebra
and $I$ an ideal of $R$ we get a map $H^i_I(R) \to H^i_I(S)$
for all $i$,  and $H^i_I(S)$ may be identified with $H^i_{IS}(S)$.
In particular, we may take $S = R$ and let the map $R \to S$
be the Frobenius endomorphism.  Since  $IS = I^{[p]}$ here,
this gives a map $H^i_I(R) \to H^i_{I^{[p]}}(R)$.  But since
$\Rad(I^{[p]}) = \Rad(I)$,  $H^i_{I^{[p]}}(R) \cong H^i_I(R)$ 
canonically. The map $H^i_I(R) \to H^i_I(R)$ so obtained again 
gives the {\it action of the Frobenius 
endomorphism} on $H^i_I(R)$.  We shall 
denote this action by $F$: note that $F(ru) = r^pF(u)$. \enddemo

\demb{Definition 2.3}  When $R$ has prime characteristic
$p>0$,  we may construct a non-commutative, associative ring $\rf$ from
$R$ which is an $R$-free left module on the symbols
$1,\,F,\,F^2, \,\ldots,\,F^e,\,\ldots$ by requiring that
$Fr = r^pF$ when $r \in R$.
We shall say that an $R$-module $M$
is an $\rf$-{\it module} if there is given an action $F:M \to M$ such
that  for all $r \in R$ and for all $u \in M$,  $F(ru) = r^pu$.
This is equivalent to the condition that $M$ be an $\rf$-module
so as to extend the  $R$-module structure on $M$. 
We then call an $R$-submodule $N$ of $M$ F-{\it stable} if
$F(N) \inc N$,  which is equivalent to requiring that $N$
be an $\rf$-submodule of $N$. If 
$M$ is any $\rf$-module and $S$ is an $R$-algebra then
there is an $S\f$-module structure on $S \otimes _R M$ determined
by the condition that $F(s \otimes u) = s^p \otimes F(u)$.

In particular, since we have an $\rf$-module structure on
$H^i_I(R)$,  we may refer to the F-stable submodules of $H^i_I(R)$.  

If $R$ is local of Krull dimension $d$ and  $\vect x d$
is a system of parameters,  then $H^d_m(R)$ may be identified
with $\dlim t R/(x_1^t,\,\ldots,\,x_d^t)$,  where the $t\,$th map
in the direct limit system 
$$ R/(x_1^t,\,\ldots,\,x_d^t)
\to R/(x_1^{t+1},\,\ldots,\,x_d^{t+1})$$
is induced by multiplication by $x_1\,\cdots\,x_d$. If $R$ is
Cohen-Macaulay the maps in this direct limit system are injective.
When $H^d_m(R)$ is thought of as a direct limit in this way,
we write $\<r;\,x_1^t, \, \ldots, \, x_d^t\>$  for the image
in $H^d_m(R)$ of the element represented by $r$ in 
$R/(x_1^t, \, \ldots, \, x_d^t)$.  The action of the Frobenius
endomorphism on the highest local cohomology module in this 
case may be described as sending 
$$\<r;\, x_1^t,\,\ldots,\,x_d^t\> \mapsto \<r^p;\, x_1^{pt},\,\ldots,\,x_d^{pt}\>.$$ 
 \enddemo

\demb{Discussion 2.4: another point of view for F-stable submodules}
Let $(R, \, m, \, K)$ be local of Krull dimension $d$, where $R$
has characteristic $p > 0$. Consider an F-stable submodule
$N \inc H^d_m(R)$.  Suppose that $R$ is reduced.  We have an isomorphism
of $(R, \, m, \, K)$ with $(S,\, n,\, L)$ where $S = R^{1/p}$ given
by $\Phi: R \to R^{1/p}$,  where $\Phi(r) = r^{1/p}$.
We have a commutative diagram:  
$$\CD  R @>{F}>> R \\
     @V{=}VV      @VV{\Phi}V\\
      R @>>{\iota}> R^{1/p} \endCD \,\, ,
$$
where $\iota: R \inc R^{1/p}$ is the inclusion map.  In general,
when $\Phi:R \to S$ is any ring isomorphism,  for each submodule
$N$ of $H^i_I(R)$  there is a corresponding submodule $N'$  of
$H^i_{\Phi(I)}(S)$.  In fact, if $\Psi=\Phi^{-1}$,
and we use $\ps Q$  to indicate restriction of scalars
from $R$-modules to $S$-modules, then $H^i_{\Phi(I)}(S)$ is
canonically isomorphic with $\ps\bigl(H^i_I(R)\bigr)$ and
$N'$ is the image of $\ps N$ in $H^i_{\Phi(I)}(S)$.  Note
that ${}^\Psi \blank$ is an exact functor. 

When $S = R^{1/p}$ and and $I=m$,  the modules  
$H^i_{\Phi(m)}(S)$, $H^i_n(S)$,  $H^i_{mS}(S)$, and $H^i_m(S)$ may all
be identified:  the
first three may be identified because $\Phi(m)$ and $mS$  both
have radical $n$, and the last two because if
$\vect f h$ generate $m$ their images $\vect g h$ in $S$ generate
$mS$ and the complexes $C^\bullet(f;\,S)$ and $C^\bullet(g;\,S)$
are isomorphic.
The condition that $N$ is F-stable is equivalent
to the condition that $N$ maps into $N'$ in $H^i_m(S) \cong H^i_{\Phi(m)}(S)$.
A very important observation is this:

\noi $(**)$ 
With notation as just above, if $N$ is F-stable and
$J \inc R$ kills $N$ (e.g., if $J = \Ann_R N$) then
$\Phi(J)$ kills the image of $N$ in $H^i_m(R^{1/p})$.

The hypothesis that $N$ is F-stable means that $N$ maps
into the corresponding submodule $N'$,  and $N'$ is
clearly killed by $\Phi(J)$.  \enddemo

\demb{Definition 2.5}  A local ring $(R,\,m)$ of
Krull dimension $d$ is 
FH-{\it finite} if, for all $i$, $0 \leq i \leq d$,  only finitely
many $R$-submodules of $H^i_m(R)$ are F-stable. We shall
that $R$ has {\it \fh} if for all $i$,  $H^i_m(R)$
has finite length in the category of $\rf$-modules. \enddemo

Our main focus in studying the properties
of being FH-finite and of having \fh\  is when the local ring $R$ 
is Cohen-Macaulay.  
Of course, in this case there is only one nonzero local
cohomology module, $H^d_m(R)$.   However, in \S5 we show
that every face ring has \fh.

Since every $H^i_m(R)$ has DCC even in the category of $R$-modules,
we know that $H^i_m(R)$ has finite length in the category of
$\rf$-modules if and only if it has ACC in the category of
$\rf$-modules.  Of course, it is also equivalent to assert that
there is a finite filtration of $H^i_m(R)$ whose 
factors are simple $\rf$-modules.  

\demb{Discussion 2.6: purity}
Recall that a map of $R$-modules $N \to N'$ is {\it pure} if for every
$R$-module $M$ the map $N \otimes_R M \to N' \otimes_R M$ is injective.
Of course, this implies that $N \to N'$ is injective, and may
be thought of as a weakening of the condition that $0 \to N \to N'$
split,  i.e., that $N$ be a direct summand of $N'$.  If $N'/N$ is
finitely presented,  $N \to N'$ is pure if and only if it is split.
For a treatment of the properties of purity, see, for example, 
[HH5, Lemma (2.1), p.\ 49].   

An $R$ algebra $S$ is called {\it pure} if $R \to S$ is pure as
a map of $R$-modules, i.e., for every $R$-module $M$,  the
map $M = R \otimes_R M \to S \otimes_R M$ is injective.  
A Noetherian ring $R$ of characteristic $p$ is
called F-{\it pure} (respectively, F-{\it split}) if the Frobenius 
endomorphism $F: R \to R$
is pure (respectively, split).  Evidently, an F-split ring is
F-pure and and an F-pure ring is reduced.  If $R$ is an F-finite
Noetherian ring,  F-pure and F-split are equivalent (since the
cokernel of $F:R \to R$ is finitely presented as a module over
the left hand copy of $R$),  and the two notions are also equivalent
when $(R,\,m,\,K)$ is complete local, for in this case, $R \to S$ is split 
iff $R \otimes_R E \to S \otimes_R E$ is injective, where
$E = E_R(K)$.  An equivalent condition is that the map obtained
by applying $\Hom_R(\blank,\,E)$ be surjective, and since
$R \cong \Hom_R(E,\,E)$, by the adjointness of tensor and
Hom that map can be identified with the maps 
$\Hom_R(S,\,R) \to \Hom_R(R,\,R) \cong R$.

We say that a local ring $R$ is F-{\it injective} if $F$ acts injectively
on all of the local cohomology modules of $R$ with support in $m$.
This holds if $R$ is F-pure. 

When $R$ is reduced, 
the map $F:R \to R$ may be identified with the algebra
inclusion $R \inc R^{1/p}$, and so $R$ is F-pure (respectively, F-split) 
if and only if it is reduced and the map $R \inc R^{1/p}$ is pure
(respectively, split). \enddemo         

\proclaim{Lemma 2.7} Let $(R,m,K)$ be a Noetherian local ring
of positive prime characteristic $p$ and Krull dimension $d$.   
{\pa
\part{a} $R$ is FH-finite (respectively, has \fh) 
if and only if its completion $\R$ is FH-finite.

\part{b} Suppose that  $(R, \,m) \to (S,\, m_S)$ is a local  
homomorphism of local rings such that
$mS$ is primary to the maximal ideal of $S$,  i.e., such
that the closed fiber $S/mS$ has Krull dimension 0.  
Suppose either that $R \to S$ is flat (hence, faithfully flat),
split over $R$, or that $S$ is pure over $R$. If $S$ is FH-finite,
then $R$ is FH-finite, and if $S$ has finite FH-length,
then  $R$ has \fh.  More generally, the poset of F-stable submodules
of any local cohomology module $H^i_m(R)$ injects in order-preserving
fashion into the poset of F-stable submodules of $H^i_{m_S}(S)$.   

\part{c} $R$ is F-injective if and
only if $\R$ is F-injective.  $R$ is F-pure if and only if
$\R$ is F-pure. \par
}
\endproclaim

\demb{Proof} Completion does not affect either what the local 
cohomology modules are nor what the action of Frobenius is.  Since
each element of a local cohomology module over $R$ is killed by
a power of $m$,  these are already $\R$-modules.  Thus, (a) is
obvious.  

Part (b) follows from the fact that the local cohomology modules
of $S$ may be obtained by applying $S \otimes_R \blank$ to those
of $R$, and that the action of $F$ is then the one discussed in 
Definition 2.3 for tensor products, i.e.,  
$F(s \otimes u) = s^p \otimes F(u)$.  From this one sees that if
$N$ is F-stable in $H^i_m(R)$,  then  $S \otimes_R N$ is F-stable
in $S \otimes_R H^i_m(R) \cong H^i_{m_S}(S)$.  Thus, we only
need to see that if $N \inc N'$ are distinct F-stable submodules
of  $H^i_m(R)$,  then the images of $S \otimes N$ and $S \otimes N'$
are distinct in $S \otimes H^i_m(R)$.  It suffices to see this
when  $R \to S$ is pure:  the hypothesis of faithful flatness or
that $R \to S$ is split over $R$ implies purity.  But  $N'/N$ injects
into  $H^i_m(R)/N$,  and $S \otimes _R N'/N$ in turn injects into 
$S \otimes_R \(H^i_m(R)/N\)$  by purity, so that the image of 
$u \in N' - N$ is nonzero in 
$$S \otimes_R \(H^i_m(R)/N\) \cong  H^i_m(S)/ \im(S \otimes N).$$  
This shows that $1 \otimes u$
is in the image of  $S \otimes N'$ in  $H^i_m(S)$ but not in
the image of  $S \otimes N$. 

Part (c), in the case of F-injectivity, follows from the fact
that it is equivalent to the injectivity of the action of F
on the  $H^i_m(R)$,  and that neither these modules nor the action of
$F$ changes when we complete.  In the case of F-purity, we prove
that if $R$ is F-pure then so is $\R$:  the other direction is trivial.  
Consider an ideal $I$ of the completion, and suppose that
there is some element $u$ of the completion such that
$u \notin I$ but $u^p \in I\p$. Choose $N$ such that $u \notin  I+m^N\R$.
We see that we may assume that $I$ is primary to the maximal ideal of
$\R$,  which implies that it is the expansion of its contraction $J$
to $R$.  Then we may choose $v \in R$ such that $v-u \in I = J\R$. 
But then $v \notin J$ but $v^p - u^p \in J\p\R$,  and since
$u^p \in J^p\R$ we have that $v^p \in J\p\R \cap R = J\p$ and
so $v \in J$,  a contradiction.  Thus, ideals of $\R$ are contracted
with respect to Frobenius, and, consequently, $\R$ is reduced.  
Then $\R \to \R^{1/p}$ is cyclically pure, by the contractedness of ideals
with respect to Frobenius that we just proved, which shows that it
is pure:  see [Ho2], Theorem (1.7).  It follows that $\R$ is F-pure.  \qed \enddemo

\demb{Discussion 2.8: the gamma construction}  Let $K$ be
a field of positive characteristic $p$ with a $p$-base $\Lambda$.
Let $\Gamma$ be a fixed cofinite subset of $\Lambda$. For $e \in \N$  
we denote by $K_{\Gamma,e}$ the purely inseparable field extension of $K$
that is the result of adjoining $p^e\,$th roots of all elements in $\Gamma$ 
to $K$,  which is unique up to unique isomorphism over $K$.  

Now suppose that $(R,\,m)$ is a complete local ring of positive prime
characteristic $p$ and that $K \inc R$ is a coefficient field,
i.e., it maps bijectively onto $R/m$.  Let $\vect x d$ be a system
of parameters for $R$,  so that $R$ is module-finite over
$A = K[[\vect x d]] \inc R$.  Let $A_\Gamma$ denote
$$\bigcup_{e \in \N} \, K_{\Gamma\!,e}[[\vect x d]],$$
which is a regular local ring that is faithfully flat and
purely inseparable over $A$.  Moreover, the maximal ideal of $A$
expands to that of $A_\Gamma$.    
We shall let $R_\Gamma$ denote $A_\Gamma \otimes _A R$,  which
is module-finite over the regular ring $A_\Gamma$ and which is
faithfully flat and purely inseparable over $R$.  The maximal
ideal of  $R$ expands to the maximal ideal of $R_\Gamma$.  The
residue class field of $R_\Gamma$ is $K_\Gamma$.  

We note that $R_\Gamma$ 
depends on the choice of coefficient field $K$ for $R$,  and the choice
of $\Gamma$,  but does not depend on the choice of system of parameters
$\vect x d$.  We refer the reader to \S 6 of [HH3] for more details.
It is of great importance that $R_\Gamma$ is F-finite, i.e., finitely
generated as a module over $F(R_\Gamma)$.  This implies that it
is excellent: see [Ku2].

It is shown in [HH3] that, if $R$ is reduced, then 
for any sufficiently small
choice of the co-finite subset $\Gamma$ of $\Lambda$, $R_\Gamma$ is
reduced.  It is also shown in [HH3] that if $R$ is
Cohen-Macaulay (respectively, Gorenstein), then $R_\Gamma$ is 
Cohen-Macaulay (respectively, Gorenstein).    \enddemo

\proclaim{Lemma 2.9} Let $R$ be a complete local ring of positive
prime characteristic $p$.  Fix a coefficient field $K$ and a
$p$-base $\Lambda$ for $K$. Let notation be as in the preceding 
discussion.  
{\pa
\part{a} Let $W$ be an Artinian $R$-module with an $\rf$-module structure
such that the action of $F$ is injective.  Then for any sufficiently
small choice of $\Gamma$ co-finite in $\Lambda$,  the action of
$F$ on $R_\Gamma \otimes_R W$ is also injective.  

\part{b} Suppose that $F$ acts injectively on a given local cohomology
module of $R$.  Then $F$ acts injectively on the corresponding local
cohomology module of $R_\Gamma$ for all sufficiently small cofinite
$\Gamma$.  In particular, if $R$ is F-injective, then so is $R_\Gamma$.

\part{c} Suppose that $R$ is F-pure.  Then for any choice of
$\Gamma$ cofinite in the $p$-base such that $R_\Gamma$ is 
reduced, and, hence, for all sufficiently small co-finite $\Gamma$,  
$R_\Gamma$ is F-pure.
\par
}
\endproclaim

\demb{Proof}  For part (a),  let $V$ denote the finite-dimensional
$K$-vector space that is the socle of $W$.  Let $W_\Gamma = R_\Gamma
\otimes_R W$.  Because the maximal ideal $m$ of $R$ expands to
the maximal ideal of $R_\Gamma$,  and  $R_\Gamma$ is $R$-flat,
the socle in $W_\Gamma$ may be identified with  
$V_\Gamma = K_\Gamma \otimes V$.  If $F$ has a nonzero kernel
on $W_\Gamma$ then that kernel has nonzero intersection with 
$V_\Gamma$, and that intersection will be some $K$-subspace of
$V_\Gamma$. Pick $\Gamma$ such that the dimension of the kernel
is minimum.  Then the kernel is a nonzero subspace  $T$ of 
$V_\Gamma$ whose intersection with $V \inc V_\Gamma$ is 0.  
Choose a basis $\vect v h$ for $V$ and choose a basis 
for $T$ as well. Write each basis vector for $T$ as $\sum_{j=1}^h a_{ij}v_j$,
where the $a_{ij}$ are elements of $K_\Gamma$.
Thus, the rows of the matrix $\alpha = (a_{ij})$ represent a basis for $T$.
Put the matrix $\alpha$ in reduced row echelon form:  the leftmost nonzero
entries of the rows are each 1,  the columns of these entries are distinct,
proceeding from left to right as the index of the row increases, and
each column containing the leading 1 of a row  has its other
entries equal to 0.  This matrix is uniquely determined by the
subspace $T$.  It has at least one coefficient $a$ not in $K$ (in fact, at
least one in every row), since $T$ does not meet $V$.  

Now choose 
$\Gamma' \inc \Gamma$ such that
$a \notin K_{\Gamma'}$,  which is possible by Lemma (6.12) of [HH3].
Then the intersection of $T$ with $V_{\Gamma'}$ must be smaller than
$T$,  or else $T$ will have a $K_\Gamma$ basis consisting of linear
combinations of the $v_j$ with coefficients in $K_{\Gamma'}$,  and this
will give a matrix $\beta$ over $K_{\Gamma'}$ with the same row space over
$K_\Gamma$ as before.  When we put $\beta$ in row
echelon form, it must agree with $\alpha$, which forces the 
conclusion that $a \in K_{\Gamma'}$, a contradiction.  

Part (b) follows  immediately from part (a).

To prove part (c),  consider a choice of  
$\Gamma$ sufficiently small that $R_\Gamma$ is reduced.
Let  $E$ be the injective hull of $K$ over $R$.  For each power
$m^t$ of the maximal ideal of $R$,  we have that $R_\Gamma/(m^t) \cong
K_\Gamma \otimes _K R/m^t$.  Thus, the injective hull of
$K_\Gamma$ over $R_\Gamma$ may be identified with 
$K_\Gamma \otimes _K E$.  We are given that the map
$E \to E \otimes_R R^{1/p}$ is injective.  We want to show
that the map 
$$
E_\Gamma \to E_\Gamma \otimes_ {R_\Gamma} (R_\Gamma)^{1/p}
$$
is injective.   Since the image of a socle generator in $E$ 
is a socle generator in  $E_\Gamma$, it is
equivalent to show the injectivity of the map
$E \to E_\Gamma \otimes_{R_\Gamma} (R_\Gamma)^{1/p}$.  

The completion of $R_\Gamma$ may be thought of as the complete tensor
product of $K_\Gamma$ with $R$ over $K \inc R$. However,  if one tensors with 
a module in which every element is killed by a power of the maximal
ideal we may substitute the ordinary tensor product for the complete tensor
product.  Moreover, since $R_\Gamma$ is reduced, 
we may identify $(R_\Gamma)^{1/p}$ with $(R^{1/p})_{\Gamma^{1/p}}$:  the
latter notation means that we are using $K^{1/p}$ as a coefficient
field for $R^{1/p}$, that we are using the $p\,$th roots $\Lambda^{1/p}$ of 
the elements of the
$p$-base $\Lambda$ (chosen for $K$) as a $p$-base for $K^{1/p}$, and that 
we are using the set $\Gamma^{1/p}$ of $p\,$th roots  of elements of 
$\Gamma$ as the co-finite subset of $\Lambda^{1/p}$ in the construction of 
$R^{1/p}_{\Gamma^{1/p}}$. But $(K^{1/p})_{\Gamma^{1/p}} 
\cong (K_\Gamma)^{1/p}$.   Keeping in mind that every element
of $E_\Gamma$ is killed by a power of the maximal ideal, and
that $E_\Gamma \cong K_\Gamma \otimes_K E$,  we have that 
$$
E_\Gamma \otimes_{R_\Gamma}R_\Gamma^{1/p} \cong 
(K_\Gamma \otimes_K E) \otimes _{K_\Gamma \otimes_K R} 
(K_\Gamma^{1/p} \otimes_{K^{1/p}} R^{1/p})
$$ 
and so, writing $L$ for $K_\Gamma$,  we have that
$$
E_\Gamma \otimes_{R_\Gamma}R_\Gamma^{1/p} \cong
(L \otimes_K E) \otimes_{L \otimes_K R} (L^{1/p} \otimes_{K^{1/p}} R^{1/p}).
$$
Now, if $K$ is any ring, $L$ and $R$ are any $K$-algebras,
$S$ is any $(L \otimes_K R)$-algebra (in our case,
$S = L^{1/p} \otimes_{K^{1/p}} R^{1/p}$), and
$E$ is any $R$-module,
there is an isomorphism 
$$
(L \otimes _K E) \otimes_{L \otimes _K R} S \cong E \otimes_R S
$$ 
which maps $(c \otimes u) \otimes s$ to $u \otimes cs$.  (The
inverse map sends $u \otimes s$ to $(1 \otimes u) \otimes s$.
Note that $(c \otimes u) \otimes s = (1 \otimes u) \otimes cs$ in
$(L \otimes _K E) \otimes_{L \otimes _K R} S$.)

Applying this fact, we find that $E_\Gamma \otimes_{R_\Gamma}R_\Gamma^{1/p}$
is isomorphic with
$$
E \otimes_R (L^{1/p} \otimes_{K^{1/p}} R^{1/p}) \cong 
E \otimes_R (R^{1/p} \otimes_{K^{1/p}} L^{1/p}) \cong
(E \otimes_R R^{1/p}) \otimes_ {K^{1/p}} L^{1/p}
$$
by the commutativity and associativity of tensor product.  But
$E$ injects into $E \otimes_R R^{1/p}$ by hypothesis and the
latter injects into  $(E \otimes_R R^{1/p}) \otimes_ {K^{1/p}} L^{1/p}$
simply because $K^{1/p}$ is a field  and $L^{1/p}$ is a nonzero
free module over it.  \qed \enddemo

\demb{Discussion 2.10:  finiteness conditions on local cohomology as
an F-module and tight closure}  We want to make some connections between 
F-submodules of local cohomology and tight closure theory.  
Let $R$ be a reduced local ring of characteristic $p > 0$.   
Let us call a submodule of $H = H^d_m(R)$ {\it strongly proper}
if it is annihilated by a nonzerodivisor of $R$.  Assume that
$R$ has test elements.  The {\it finitistic} tight closure of
$0$ in a module $M$ is the union of the submodules  $0^*_N$  as
$N$ runs through the finitely generated submodules of $M$. 
It is not known, in general, whether the tight closure of 0 in an Artinian
module over a complete local ring is the same as the finitistic
tight closure:  a priori, it might be larger.  Cf.\ [LySm1,2], [El], and 
[St] for results in this direction. 

However, if $\rmk$
is an excellent reduced equidimensional local ring with $\di(R) = d$,  the two are
the same for $H^d_m(R)$:  if  $u$ is in the tight closure of  0  and represented
by  $f$  mod  $I_t = (x_1^t, \, \ldots,  x_d^t)$ in  $H^d_m(R)  = 
\dlim t R/(x_1^t, \, \ldots, \, x_d^t)R$,  then there exists
$c \in R\0$ such that for all $q = p^e \gg 1$ the class of $cu^q$  
maps to 0  under the map $R/I_{qt}  
\to H^d_m(R)$, i.e., for some  $k_q$,  $cu^q(x_1\, \cdots x_d)^{k_q}
\in I_{qt+k+q}$ for all $q \gg 1$.
But $I_{qt+k_q}:_R  (x_1\, \cdots x_d)^{k_q} \inc I_{qt}^*$
(since $R$ is excellent, reduced, and local, it has a completely stable
test element, and this reduces to the complete case, which follows
from [HH2], Theorem (7.15)),  and so if $d$ is a test element for  $R$  we have
that  $cdu^q \in I_{qt} = I_t^{[q]}$ for all $q \gg 0$, and
so the class of $u$ mod $I_t$  is in the tight closure of
0 in  $R/I_t$  and hence in the image of $R/I_t$ in $H^d_m(R)$,  
as required.

Let us note that the finitistic tight closure of 0 in $H$ is an
an F-stable strongly proper submodule of $H^d_m(R)$, as shown in [LySm2], Proposition (4.2).  The
reason is that it is immediate from the definition of
tight closure that if $u \in 0^*_N$  then  $u^q \in 0^*_{N^{[q]}}$,
where  $q = p^e$ and $N^{[q]}$ denotes an image of $F^e(N)$
in  $F^e(H)$  for some ambient module  $H \supseteq N$.  In
particular,  $u^q \in 0^*_{F^e(N)}$.  Moreover,  if $c$ is
a test element for the reduced ring $R$,  then $c \in R\0$ and
so $c$ is a nonzerodivisor,  and  $c$ kills  $0^*_N$  for
every finitely generated $R$-module $N$.    

Conversely,
any strongly proper F-stable submodule $N \inc H$ is in the tight closure
of 0.  If $c$ is a nonzerodivisor that kills  $N$  and
$u \in N$,  then $cu^q = 0$  for all $q$:  when we identify
$F^e(H)$  with $H$,   $u^q$ is identified with $F^e(u)$.  

What are the strongly proper submodules of $H = H^d_m(R)$?  
If $(R,\,m,\,K)$ is complete with  $E = E_R(K)$ the injective hull of
the residue class field and 
canonical module $\w := \Hom_R(H^d_m(R), \, E)$,   
then submodules of $H$ correspond to 
the proper homomorphic images of $\w$:  the inclusion
$N \inc H$  is dual under $\Hom_R(\blank,\,E)$ to a surjection  
$\w \surj \Hom_R(N,\,E)$.
If $R$ is a domain, for every proper $N \inc H$ we have
that  $\w \surj \Hom_R(N,\,E)$ is a proper surjection,
and therefore is killed by a nonzerodivisor. Therefore, we have
the following results of K.\ E.\ Smith [Sm3]  (see Proposition (2.5),
p.\ 169 and the remark on p.\ 170 immediately following the
proof of (2.5)). (See also [Sm1], Theorem (3.1.4), where the restricted
generality is not needed.)

\proclaim{Proposition 2.11 (K.\ E.\ Smith)}  If  $R$ is a reduced 
equidimensional excellent local ring of characteristic $p$,  
then the tight closure $0^*$
of $0$ in $H^d_m(R)$  (which is the same in the finitistic and
ordinary senses) is the largest strongly proper F-stable submodule of
$H^d_m(R)$. If $R$ is a complete local domain, it is the largest
proper F-stable submodule of $H^d_m(R)$. \qed \endproclaim  
\bigskip

A Noetherian local ring is called F-{\it rational} if some
(equivalently, every) ideal generated by parameters is tightly closed.
An excellent F-rational local ring is a Cohen-Macaulay 
normal domain. The completion of an excellent F-rational local
ring is again F-rational. Cf.\ [HH3], Proposition (6.27a).   
From this and the discussion 
above we have at once (see [Sm3], Theorem (2.6), p.\ 170):  

\proclaim{Proposition 2.12 (K.\ E.\ Smith)} Let $R$ be an excellent
Cohen-Macaulay local ring of characteristic $p$ and Krull dimension
$d$.  Then $R$ is F-rational if and only if $H^d_m(R)$ is a simple
$\rf$-module. \qed \endproclaim

\demb{Example 2.13}  The ring obtained by killing the size $t$ minors
of a matrix of indeterminates in the polynomial ring in those 
indeterminates is an example of an F-rational ring. In fact this ring is weakly F-regular, i.e.,
every ideal is tightly closed.  The local ring at the origin is
therefore FH-finite by the above result:  the unique non-vanishing local cohomology 
module is $\rf$-simple. \enddemo

\proclaim{Proposition 2.14}(see also (4.17.1) in [Sm4]). Let $R$ be a Cohen-Macaulay local
domain and suppose that there is an $m$-primary ideal $\fA$
such that   $\fA I^* \inc I$ for every ideal $I$ of $R$
generated by part of a system of parameters.  Then $R$ has \fh.\endproclaim
\demb{Proof} Let $d$ be the dimension of $R$.  By the discussion above,
every proper F-stable submodule of $H = H^d_m(R)$ is contained
in $0^*_H$.  But the discussion above shows that $0^*_H$ is
a union of submodules of the form $I^*/I$ where $I$ is a parameter
ideal, and so $0^*_H$ is killed by $\fA$,  and has finite length 
even as an $R$-module.  \qed 

See Theorem (4.21),  which gives a stronger conclusion when the residue
class field is perfect and $R$ is F-injective.

\demb{Example 2.15} Let $R = K[[X,\,Y,\,Z]]/(X^3 + Y^3+ Z^3)$, 
where $K$ is a field of positive characteristic different from
3.  Then $R$ is a Gorenstein domain,  and the tight closure of 0  in
$H^2_m(R)$ is just the socle, a copy of $K$:  the tight closure
of every parameter ideal is known to contain just one additional
element, a representative of the generator of the socle modulo
the parameter ideal.  Evidently $R$ is FH-finite.  It is known 
(see, for example, [HR], Prop. (5.21c), p.\ 157)
that $R$ is F-injective if and only if the characteristic of
$K$ is congruent to 1 modulo 3.    If the characteristic
is congruent to 2 modulo 3,  $R[[t]]$ does not have FH-finite
length by Theorem (4.15) of \S4.

%%  -on page 13, line 17, K=k(u,v).   

\demb{Example 2.16}  We construct a complete local F-injective domain
of dimension one (hence, it is Cohen-Macaulay) that is not FH-finite.  Note
that Theorem (4.21) implies that there are no such examples when
the residue class field of the ring is perfect.

Let $K$ be an infinite field of characteristic $p > 0$ (it will be necessary
that $K$ not be perfect)  and let $L$ be a finite 
algebraic extension field of $K$   such that
{\pa  
\part 1 $[L : L^p[K]] > 2$   (all one needs is that  the dimension of
    of   $L/L^p[K]$   over  $K$ is at least  2) and\medskip  

\part 2 $L$  does not contain any element of  $K^{1/p} - K$
    (equivalently,   $L^p  \cap  K = K^p$). \par}
    \medskip  

Then the quotient $L/K$  has infinitely many $K\f$-submodules but $F$ acts 
injectively on it.  Moreover,  if  $R = K + xL[[x]] \inc L[[x]]$  then
$R$  is a complete local  one-dimensional domain that is F-injective but 
not FH-finite.\medskip  

The conditions in (1) and (2) above may be satisfied as follows: 
if  $k$  is infinite perfect of characteristic $p > 2$,  $K = k(u,v)$,
where $u$ and $v$ are indeterminates, and
$$L = K[y]/(y^{2p} + uy^p - v),$$ then   (1) and (2) above are satisfied.\medskip  

\bpf  The image of   $L$  under  $F$  is  $L^p$ --- this need not be
a  $K$-vector space,  but  $L_1 =  L^p[K]$   is a  $K$-vector space 
containing the image of  $F$.  All of the  $K$-vector subspaces 
of  $L$  strictly
between  $L_1$   and  $L$   are  F-stable,  and there are infinitely
many.   The statement that  $F$  acts injectively on  $L/K$  is exactly
the statement that  $L^p \cap K = K^p$.  \medskip  

With  $R$  as above, the exact sequence   
$$0  \to R \to L[[x]] \to L/K \to 0$$ yields a long exact sequence for
local cohomology:
$$0 \to L/K \to H^1_m(R) \to H^1_m(L[[x]]) \to 0,$$
Since  $L/K$  embeds in  $H^1_m(R)$  as an  F-stable submodule
and  $m$   kills it,  its  $R$-module structure is given by
its  $K$-vector space structure.  Moreover, since  $F$  is injective
on  $H^1_m(L[[x]])$,  F-injectivity holds for  $R$  iff it holds
for  $L/K$.  \medskip   

This establishes all assertions except that the given example 
satisfies (1) and (2).  Note that the equation  $y^{2p} + uy^p - v$
is irreducible over  $k[y,u,v]$  (the quotient is  $k[y,u]$).  
Suppose that  $L$  contains an element  $w$  of  $K^{1/p}$  not in  $K$.
Then  $[K[w] :  K] = p$,   and so   $[L : K[w]] = 2$.  It follows
that $y$  satisfies a monic quadratic equation over  $K[w]$.  \medskip  

But if we enlarge  $K[w]$  to all of  $K^{1/p}$  we know the quadratic
equation that  $y$  satisfies:   $y^2 + u^{1/p}y - v^{1/p} = 0$,
which is clearly irreducible over $K^{1/p} = k(u^{1/p}, v^{1/p})$.  
This quadratic is unique, so we must have  $u^{1/p},\,  v^{1/p}$  are both
in  $K[w]$,  a contradiction,  since  adjoining both produces an
extension of  $K$  of degree  $p^2$.  \medskip  

It remains to determine $[L : L^p[K] ] = [K[y] : K[y^p]]$.
Since  $y^p$  satisfies an irreducible quadratic equation over  $K$,
$[K[y^p] : K] = 2$,  and so  $[K[y] : K[y^p] ] = 2p/2 = p > 2$,  by
assumption. \qed

\quad\bigskip
\bigskip
{\chbf 
\centerline{\hdbf 3. ANNIHILATORS OF F-STABLE SUBMODULES AND THE} 
\smallskip
\centerline{\hdbf FH-FINITE PROPERTY FOR F-PURE GORENSTEIN LOCAL RINGS} 
}
\bigskip

In this section we shall prove a theorem of independent interest which can be used to establish that certain families of
radical ideals in excellent local rings are finite. As an immediate corollary to it, we obtain that
if $R$ is local, excellent, then any family of radical ideals closed under sum, intersection and
primary decomposition is finite. Another consequence is that if $R$ is a local ring
of positive prime characteristic $p$ and  $M$ is an Artinian
$\rf$-module such that $F$ acts injectively on $M$, then
the set of annihilator ideals in $R$ of F-stable submodules
of $M$ is a finite set of radical ideals closed under primary
decomposition (R.~Y.~Sharp proved this result independently: see Theorem (3.10) and Corollary (3.11) in [Sh]).  In fact, it consists of a finite set of prime
ideals and their intersections.  From this we deduce that an F-pure Gorenstein local ring is FH-finite.
We say that a family of radical ideals of a Noetherian ring
is {\it closed under primary decomposition} if for every
ideal $I$ in the family and every minimal prime $P$ of $I$,
the ideal $P$ is also in the family. 

The following result is the main Theorem of this section. 

\proclaim{Theorem 3.1} 
Let  $M$  be a Noetherian module over an excellent local ring  $(R,\,m)$.
Then there is no family $\{N_\lambda\}_{\lambda \in \Lambda}$ of submodules 
of $M$ 
satisfying all four of the conditions below: \smallskip
{\pa  
\part 1 the family is closed under finite sums\smallskip
\part 2 the family is closed under finite intersection\smallskip
\part 3 all of the ideals $\Ann_R(M/N)$ for $N$ in the family are
radical, and\smallskip
\part 4 there exist
   infinitely many modules in the family such that if  $N$, $N'$  are any
   two of them, the minimal primes of  $N$  are mutually incomparable
   with the minimal primes of  $N'$.  \par
   }
   \smallskip
Hence, if a family of submodules $\{N_\lambda\}_{\lambda \in \Lambda}$ 
of $M$ satisfies conditions (1), (2), and (3) above and the set 
$$\{ \Ann_R(M/N_\lambda): \lambda \in \Lambda\}$$ 
is closed under primary 
decomposition, then this set of annihilators is finite. 
\endproclaim 
\demb{Proof}  Assume that one has a counterexample.
We use both induction on the dimension of $R$ and Noetherian induction
on $M$.    Take a counterexample in which the ring
has minimum dimension.  One can pass to the completion.
Radical ideals stay radical, and (4) is preserved (although
there may be more minimal primes).  The key point is that
if  $P$, $Q$  are incomparable primes of  $R$,  and
$\R$ is the completion of  $R$,  then  $P\R$,  $Q\R$  are radical with no
minimal prime in common.  A common minimal prime would contain
$(P + Q)\R$, a contradiction, since the minimal primes of  $P\R$
lie over $P$). This, applied together with the fact that
the minimal primes of $\widehat{M/N}$ are minimal over $P\R$ for
some minimal prime $P$ of $M/N$ enables us to pass to the completion.    

Take infinitely many  $N_i$   as in  (4).
Let  $M_0$  be maximal among submodules
of $M$ contained in infinitely many of the  $N_i$.
Then the set  of modules in the family containing  $M_0$
gives a new counterexample, and we may pass to all quotients
by  $M_0$  ($M_0$ need not be in the family to make this reduction.)
Thus, by Noetherian induction on  $M$  we may assume that every
infinite subset of the  $N_i$  has intersection 0.  

Consider the set of all primes of  $R$  in the support of
an $M/N_i$.  If  $Q \not= m$  is in the support of infinitely many
we get a new counterexample over  $R_Q$. The  $(N_i)_Q$
continue to have the property that no two have a minimal 
prime in common (in particular, they are distinct).  Since  $R$  had 
minimum dimension for a counterexample,  we can conclude that every 
$Q$ other than   $m$  is in the support of just finitely many  $M/N_i$.

Choose  $h$   as large as possible such that there infinitely
many primes of height  $h$   occurring among the minimal primes
of an  $M/N_i$.   Then there are only finitely many primes of 
height $h+1$  or more occurring as a minimal prime of an
$M/N_i$,  and, by the preceding paragraph, each one occurs
for only finitely many  $N_i$.  Delete sufficiently many  $N_i$
from the sequence so that the no prime of height bigger than
$h$  occurs among the minimal primes of the  $M/N_i$.   

Let $D_1(i) = \bigcap_{s=1}^i\,N_s$. By
Chevalley's lemma, $D_1(i_1)$  is contained in  $m^2M$  for  $i_1$
sufficiently large:  fix such a value of $i_1$.  
Let  $W_1 = D_1(i_1)$.  Let  
$ D_2(j) = \bigcap_{s=i_1+1}^j\,N_s$.
Then $D_2(j)$   is contained in  $m^2M$   for sufficiently large
$j$:  fix such a value  $i_2$.   
Recursively, we can choose a strictly increasing sequence of integers
$\{i_t\}_t$ with $i_0 = 0$ such that every 
$$W_t = \bigcap_{s=i_{t-1}+1}^{i_t}\,N_s$$ is contained in $m^2M$.     
In this way we can construct a sequence   $W_1,\, W_2,\, W_3,\,  \ldots$  
with the same properties as the $N_i$ but such that  all of them are in   
$m^2M$.   
Now $W_1 + W_2 + \, \cdots\,+ W_t$   stabilizes for  $t \gg 0$, since
$M$ is Noetherian,  and the stable value  $W$  is contained in  $m^2M$.  
There cannot be any prime other than   $m$   in the support of   
$M/W$,  or it will be in the support of  $M/W_j$  for all  $j$  and 
this will put it in the support of infinitely many of the original  
$N_j$.  Hence, the annihilator of  $M/W$  is an  $m$-primary
ideal,  and, by construction, it is contained in   $m^2$  and,
therefore, not radical, a contradiction. 

It remains only to prove the final statement. If the set of
annihilators were infinite it would contain infinitely many
prime ideals.  Since there are only finitely many possibilities
for the height, infinitely
many of them would be prime ideals of the same height. The modules
in the family having these primes as annihilator satisfy (4),
a contradiction.  
\qed \enddemo

By applying Theorem (3.1) to the family of ideals of $R$, we immediately have: 

\proclaim{Corollary 3.2}   A family of radical
ideals in an excellent local ring closed under sum, intersection, 
and primary decomposition is finite. \qed\endproclaim

\demb{Discussion 3.3}
For any local ring $(R,\,m,\,K)$ we let  $E$ denote an injective
hull of the residue class field, and we write $\blank^\lor$ for 
$\Hom_R(\blank,\,E)$.  Note that $E$ is also a choice for
$E_{\R}(K)$, and that its submodules over $R$ are the same
as its submodules over $\R$.  $E$ is determined up to non-unique
isomorphism: the obvious map $\R \to \Hom_R(E,\,E)$ is
an isomorphism, and so every automorphism of $E$ is given by
multiplication by a unit of $\R$.  

Now suppose that $R$ is complete.  
Then $R^\lor \cong E$ and $E^\lor \cong R$, by 
Matlis duality.  Matlis duality gives an anti-equivalence between modules
with ACC and modules with DCC:  in both cases, the functor used
is a restriction of $\blank^\lor$.  In particular, the natural
map $N \to N^{\lor\lor}$ is an isomorphism whenever $N$ has DCC or ACC.   
Note that there is an order-reversing bijection between ideals $I$ of $R$ 
and submodules $N$ of $E$ given by $I \mapsto \Ann_E I$  and
$N \mapsto  \Ann_R N$:  this is a consequence of the fact that
the inclusion  $N \inj E$  is dual to a surjection $N^\lor \surjb R$
so that  $N^\lor \cong R/I$  for a unique ideal $I$ of $R$,  and
since $I = \Ann_RN^\lor$,  we have that $I = \Ann_RN$. Note that
$N \cong N^{\lor\lor} \cong (R/I)^\lor \cong \Hom_R(R/I,\,E) 
\cong \Ann_EI$.  When $R$ is regular or even if $R$ is Gorenstein,  
$E \cong H^d_m(R)$. 

When $R$ is complete local and $W$ is Artinian,
Matlis duality, provides a bijection between
the submodules of $W$ and the surjections from $W^\lor = M$,
and each such surjection is determined by its kernel $N$.
This gives an order-reversing bijection between the submodules
of $W$ and the submodules of $M$.   Specifically,  $V \inc W$
corresponds to $\Ker(W^\lor \surj V^\lor) = \Ker(M \surj V^\lor)$,
and $N \inc M$ corresponds to $\Ker\bigl(M^\lor \surj (M/N)^\lor\bigr)$.
Here, $M^\lor = (W^\lor)^\lor \cong W$ canonically.  This
bijection converts sums  to intersections and intersections to
sums:  the point is that the sum (intersection) 
of  a family of submodules is the smallest (respectively, largest)
submodule containing (respectively, contained in) all of them,
and the result follows from the fact that the correspondence
is an order anti-isomorphism.  Since the annihilator of a module
kills the annihilator of its dual, Matlis duality preserves
annihilators:  it is obvious that the annihilator of a module
kills its dual, and we have that each of the two modules is
the dual of the other.  In particular, under the order-reversing
bijection between submodules $V$ of $W$ and submodules $N$ of $M$,
we have that $\Ann_R V = \Ann_R(M/N)$.  
\enddemo

\demb{Discussion 3.4} Let $(R,\,m,\,K) \to (S,\,n,\,L)$ be local,
and suppose that $mS$ is $n$-primary and that $L$ is finite
algebraic over $K$:  both these conditions hold if $S$
is module-finite over $R$.   Let $E= E_R(K)$ and
$E_S(L)$ denote choices of injective hulls for $K$ over $R$ and for 
$L$ over $S$, respectively.  The functor
$\Hom_R(\blank,\, E)$ from $S$-modules to $S$-modules may be
identified with $\Hom_R(\blank \otimes_S S, \, E) \cong
\Hom_S\bigl(\blank, \, \Hom_R(S, \,E)\bigr)$, which shows
that $\Hom_R(S, \,E)$ is injective as an $S$-module.  Every
element is killed by a power of the maximal ideal of $S$,
since $mS$ is primary to $n$,  and the value of the functor
on $L = S/n$ is $\Hom_R(L,\, E) \cong \Hom_R(L,\,K)$ 
since the image of $L$ is killed by $m$.  But this is $L$  as
an $L$-module.  Thus, $E_S(L) \cong \Hom_R(S, \, E)$,  and the
functor $\Hom_R(\blank, E)$,  on $S$-modules, is isomorphic
with the functor $\Hom_S(\blank, \, E_S(L))$.  
\enddemo

The following Proposition can be seen as a 
consequence of the more general Theorem (3.6) and its corollary in [Sh]. However, its proof is not very difficult and we include it here for the convenience of the reader.

\proclaim{Proposition 3.5} Let $R$ be a ring of characteristic
$p$ and let $W$ be an $\rf$-module.
{\pa
\part a If  $F$ acts injectively on $W$,  the annihilator
in $R$ of every F-stable submodule is radical.
\part b If $I$ is the annihilator of an F-stable submodule
$V$ of $W$,  then $I:_Rf$  is also the annihilator of
an F-stable submodule, namely,  $fV$.  Hence, if $I$ is
radical with minimal primes $\vect P k$ then every $P_j$ (and
every finite intersection of a subset of the $P_j$
is the annihilator of an F-stable submodule of $M$.\par
} \endproclaim
\demb{Proof}  If $V$ is  F-stable and $u \in R$ is such
that $u^p \in \Ann_R V$,  then $F(uV) = u^pF(V) \inc u^p V = 0$.
Since  $F$ is injectively on $W$, $uV = 0$.  This proves part
(a).  For part (b), note that  $fV$ is F-stable since
$F(fV) = f^pF(V) \inc f^pV \inc fV$  and  $u(fV) = 0$ iff
$(uf)V = 0$ iff $uf \in \Ann_R V = I$ iff $u \in I:_R f$.
For the final statement, choose $f$ in all of the $P_j$ except $P_i$,
and note that $I:_R f = P_i$.  More generally, given a
subset of the $P_j$,  choose $f$ in all of the minimal
primes except those in the specified subset. \qed\enddemo

Now we are in position to state an important consequence of our main result in this section. This result has also been obtained by Sharp in [Sh], Theorem
(3.10) and more precisely, Corollary (3.11). Our proof is via Theorem (3.1),
so we will include here.

\proclaim{Theorem 3.6 (R.\ Y.\ Sharp)} Let $R$ be a local ring
of positive prime characteristic $p$ and let $W$ be an
Artinian $\rf$-module.  Suppose that $F$ acts injectively
on $W$.  Then the 
$$\{\Ann_R V: V  \hbox{\ is an F-stable submodule of\ }W\}$$
is a finite set of radical ideals, and consists of all
intersections of the finitely many prime ideals in it.
\endproclaim
\demb{Proof}  
By Proposition (3.5),  it suffices to prove that
family of annihilators is finite.   We may
replace $R$ by its completion without changing $M$ or
the action of $F$ on $M$. The set of F-stable submodules
is unaffected.  The annihilator of each such submodule
in $R$ is obtained from its annihilator in $\R$ by intersection
with $R$.  Therefore, it suffices to prove the result when
$R$ is complete, and we henceforth assume that $R$ is complete.   

As in (3.3) fix an injective hull $E$ of
$K$ and let $\blank^\lor = \Hom_R(\blank,\,E)$.  
Matlis duality gives a bijection of submodules of $W$
with submodules of $M = W^\lor$.  The F-stable submodules
of $W$ are obviously closed under sum and intersection.
Therefore, the submodules $N$ of $M$ that correspond to
them are also closed under sum and intersection.  We
refer to these as the {\it co-stable} submodules of $M$.
The annihilators of the modules $M/N$,  where $N$ runs through
the co-stable submodules of $M$, are the same as the annihilators
of the F-stable submodules of $W$.  We may now apply the 
final statement of Theorem (3.1).
\qed\enddemo

It is now easy to prove the second main result of this section.  Recall
the a local ring $\rmk$ of Krull dimension $d$ is {\it quasi-Gorenstein}
if $H^d_m(R)$ is an injective hull of $K$:  equivalently, this means
that $R$ is a canonical module for $R$ in the sense that its Matlis
dual is $H^d_m(R)$.

\proclaim{Theorem 3.7} Let $(R, \, m,\, K)$ be a local ring 
of prime characteristic $p > 0$.  If $R$ is F-pure
and quasi-Gorenstein, then $H^d_m(R)$ has only finitely
many F-stable submodules.  Hence, if $R$ is F-pure and 
Gorenstein, then $R$ is FH-finite. \endproclaim
\demb{Proof}  There is no loss of generality in
replacing $R$ by its completion.
We apply Theorem (3.6) to the action of F
on $H^d_m(R) = E$.  The point is that because $E^\lor = R$,
the dual of the F-stable module $V$ has the form $R/I$,  where
$I$ is the annihilator of $V$,  and so $V$ is uniquely determined
by its annihilator.  Since there are only finitely many
possible annihilators, there are only finitely many
F-stable submodules of  $H^d_m(R)$.  \qed\enddemo

In relation to this Theorem the following observation is interesting
\footnote{We thank Karl Schwede who suggested we incorporate this remark in our paper.}. Since we do not know a reference for it in the case of F-rational rings, we will
also sketch a proof of it. For the definition of tight closure in the case of modules and its variant of finitistic tight closure, we refer 
the reader to [HH2]. The result on F-rational rings is not actually related to rest of the paper, but we 
include it for the sake of completness.

\proclaim{Remark  3.8} Let $(R, \, m,\, K)$ be a local ring 
of prime characteristic $p > 0$. Assume that $R$ is quasi-Gorenstein. Then
$R$ is F-pure if and only if $R$ is F-injective.  If  $R$ is excellent as well,
then $R$ is weakly F-regular if and only if $R$ is F-rational.  
\endproclaim
\demb{Proof}
If $R$ is F-pure (respectively, weakly F-regular), then it is immediate that $R$ is F-injective (respectively F-rational).

Assume that $R$ is F-injective. To prove that $R$ is F-pure, one can proceed exactly as in Lemma (3.3) in [Fe]. 

Now assume that $R$
is $F$-rational. We can assume that $R$ is complete. Let $E = E_R(K)$. To prove that $R$ is weakly F-regular we need to show that that the finitistic
tight closure of zero in $E$ equals zero, that is $0=0^{*,fg}_{E}$.
But $E$ is isomorphic to $H=H^d_m(R)$, and $0^{*,fg}_H = 0^*_{H} =0$ since $R$ is F-rational. This finishes the sketch of the proof.\qed\enddemo

\quad\bigskip\bigskip
{\chbf 
\centerline{\hdbf 4. F-PURITY, FINITE LENGTH, 
AND ANTI-NILPOTENT MODULES}
}
\bigskip
In this section we prove that certain quotients by annihilators
are F-split, and we study the family of F-stable submodules of
the highest local cohomology both in the F-pure Cohen-Macaulay
case, and under less restrictive hypotheses.  We do not know
an example of an F-injective ring which does
not have finite FH-length, but we have not been able to prove
that one has finite FH-length even in the F-split Cohen-Macaulay
case.  We also give various characterizations of when a
local ring has finite FH-length.

\proclaim{Theorem 4.1} Let $(R, \, m, \, K)$ be a 
local ring of prime characteristic $p > 0$ of Krull dimension $d$.
Suppose that $R$ is F-split.  
Let $N$ be an F-stable submodule
of $H^d_m(R)$, and let $J = \Ann_R \,N$.  Then $R/J$ is
F-split.  In fact, let $\Phi : R \to R^{1/p}$ such that $ \Phi(r ) = r^{1/p}$. If  $T:R^{1/p} \to R$ is any $R$-linear
splitting,  then for every such annihilator ideal $J$,  
 $T\bigl(\Phi(J)\bigr) \inc J$,   and so
$T$ induces a splitting  $(R/J)^{1/p} \cong R^{1/p}/\Phi(J)
\surj R/J$.   \endproclaim

\demb{Proof} Let $H = H^d_m(R)$.  When we apply $\blank \otimes_R H$
to  $\iota: R \inc R^{1/p}$ and to $T:R^{1/p} \to R$,  
we get maps $\alpha: H \to R^{1/p}\otimes H$ $\cong H^d_m(R^{1/p}) \cong
H^d_n(R^{1/p})$,  where $n$ is the maximal ideal of $R^{1/p}$,
and also a map $\tilde T: H^d_m(R^{1/p}) \to H$.   Therefore $\alpha$ equals $\iota \otimes_R \id_H$
and $\tilde T$ equals $T \otimes_R \id_H$.    

Let $u \in H$ and $s \in R^{1/p}$. We have $\tilde T(s\alpha(u)) = 
(T \otimes_R \id_H)\bigl(s(1 \otimes u)\bigr) = 
(T \otimes_R \id_H)(s \otimes u) = T(s) \otimes u$,  and since
$T(s) \in R$,  this is simply $T(s)u$.   
To show that $T(J^{1/p}) \inc J$,  we need to prove that
$T(J^{1/p})$ kills $N$ in $H^d_m(R)$.  Take $u \in N$ and 
$j \in J$.  Then $T(j^{1/p})u = \tilde T(j^{1/p}) (\alpha(u))$,
taking $s = j^{1/p}$.    
But now, since $N$ is an F-stable submodule of $H^d_m(R)$,
$\alpha$  maps $N$ into the corresponding submodule $N'$ of
$H^d_m(R^{1/p})$, whose annihilator in $R^{1/p}$ is
$\Phi(J)$.   
We therefore  have that $j^{1/p}$ kills $\alpha(N) \inc N'$.
This is the displayed fact $(**)$ in Discussion (2.4).
Therefore,  $T(j^{1/p})u = 0$, and $T(j^{1/p}) \in J$.  \qed\enddemo

We now want to discuss the condition of having finite FH-length.

\proclaim{Proposition 4.2} Let $R$ be a characteristic $p$ local
ring with nilradical $J$, and let $M$ be an Artinian $\rf$-module.
Then $M$ has finite $\rf$-length if and only if
$JM$ has finite length as an $R$-module and $M/JM$ has finite
$(R/J)\f$-length.  \endproclaim

\demb{Proof} $JM$ has a finite filtration by submodules $J^tM$,
and $F$ acts trivially on each factor.  \qed\enddemo

Because of Proposition (4.2), we shall mostly limit our
discussion of finite FH-length to the case where $R$ is
reduced.  

\demb{Discussion 4.3}  Let $(R,\,m,\,K)$ be a local ring
of characteristic $p > 0$ and let $M$ be an Artinian $\rf$-module.  
We note that $M$ is also an Artinian $\wh{R}\f$-module with the
same action of F,  and we henceforth assume that $R$ is complete
in this discussion.  We shall also assume that  $R$
is reduced.  (In the excellent case, completing will not
affect whether the ring is reduced.) Fix an injective hull
$E = E_R(K)$ for the residue field and let $\blank^\lor$ denote the
functor $\Hom_R (\blank, \, E)$.  \enddemo

\proclaim{Lemma 4.3} Let $(R,\,m,\,K)$ be a local ring of
characteristic $p$. Then we may construct $S$ local and faithfully flat
over $R$ with maximal ideal $mS$  such that $S$ is complete and
faithfully flat over $R$,  such that $S$ is Cohen-Macaulay if
$R$ is Cohen-Macaulay, such that $S$ is F-injective if $R$ is,
and such that $S$ is F-split if $R$ is F-pure. For every $i$,
the poset of F-stable modules of $H^i_m(R)$ injects by a strictly
order-preserving map into the poset of F-stable modules of
$H^i_{mS}(S) = H^i_m(S)$.  Hence, $R$ is FH-finite (respectively,
has finite FH-length) if $S$ is FH-finite (respectively, has FH-finite
length).  \endproclaim

\demb{Proof} By Lemma (2.7), we may first replace $R$ by its completion 
$\wh{R}$. We then  choose
a coefficient field $K$ and a $p$-base $\Lambda$ for $K$, and
replace $\wh{R}$ by $\wh{R}_\Gamma$ for $\Gamma$ a sufficiently
small cofinite subset of $\Lambda$, using Lemma (2.9).  Finally,
we replace  $\wh{R}_\Gamma$ by its completion. The map on posets
is induced by applying $S \otimes_R \blank$.\qed

\demb{Discussion 4.4: reductions in the Cohen-Macaulay F-pure case}  

Consider the following three hypotheses on a local ring $R$ of
prime characteristic $p >0$.  \smallskip
{\pa
\part 1 $R$ is F-injective and Cohen-Macaulay, with perfect residue class field.
\part 2 $R$ is F-pure.
\part 3 $R$ is F-pure and Cohen-Macaulay. \par
}\smallskip
In all three cases we do not know, for example, whether the top local
cohomology module has only finitely many F-stable submodules.   
The point we want to make is that Lemma (4.3)  permits us to reduce 
each of these questions to the case where $R$ is complete and F-finite.
Moreover, because F-pure then implies F-split, in cases (2) and (3)
the hypothesis that $R$ be F-pure may be replaced by the hypothesis
that $R$ be F-split.   
\enddemo

If $V \inc V' \inc W$  are F-stable $R$-submodules of the
$\rf$-module $W$,  we refer to  $V'/V$ as a {\it subquotient} of $W$. Also,
we let $\widetilde{F}^k(V)$ denote the $R$-span of $F^k(V)$ in $W$.
 
We next note the following:

\proclaim{Proposition 4.5}  Let $R$ be a ring of positive
prime  characteristic $p$ and let  $W$ be an $\rf$-module.
The following conditions on $W$ are equivalent:
{\pa
\part a  If $V$ is an F-stable submodule of $W$,  then
$F$ acts injectively on $W/V$
\part b  $F$ act injectively on every subquotient of $W$.
\part c The action of $F$ on any subquotient of $W$
is not nilpotent.
\part d The action of $F$ on any nonzero subquotient of $W$ is
not zero.
\part e If $V \inc V'$ are F-stable submodules of $W$ such
that $F^k(V') \subseteq V$ for some $k \geq 1$ then $V' = V$.\par } 
\endproclaim
\demb{Proof}  (a) and (b) are equivalent because
a subquotient $V'/V$ is an $\rf$-submodule of $W/V$.  If
the action of $F$ on  a subquotient $V''/V \inc W/V$ is
not injective, the kernel has the form  $V'/V$ where
$V'$ is F-stable.  Hence (b) and (d) are equivalent.
If  $F$ is nilpotent on $V''/V$ it has a nonzero kernel
of the form $V'/V$.  This shows that (c) is also equivalent.
(e) follows because $F^k$ kills $V'/V$ if and only if
$F^k(V') \inc V$.  \qed\enddemo

\demb{Definition 4.6} With $R$ and $W$ as in Proposition (4.5)
we shall say that $W$ is {\it \an} if it satisfies
the equivalent conditions (a) through (e).

The following is Theorem (4.7) on p.\ 108 of [Ly]:

\proclaim{Theorem 4.7 (Lyubeznik)} Let $R$ be a local ring of
prime characteristic $p > 0$ and let $W$ be an Artinian
$R$-module that has an action $F$ of Frobenius on it. Then
$W$ has a finite filtration  
$$
0 = L_0 \inc N_0 \inc L_1 \inc N_1 \inc \cdots \inc L_s \inc N_s=M
$$
by F-stable submodules such that every $N_j/L_j$ is nilpotent, i.e.,
killed by a single power of $F$,  while every $L_j/N_{j-1}$ is
simple in the category of $\rf$-modules, with a nonzero action
of $F$ on it.  The integer $s$ and the isomorphism classes
of the modules $L_j/N_{j-1}$  are invariants of $W$.  
\qed
\endproclaim

Note that the assumption that the action of $F$ on a simple
module $L \not=0$ is nonzero is equivalent to the
assertion that the action of $F$ is injective, for if $F$
has a non-trivial kernel it is an $\rf$-submodule and so
must be all of $L$.  

We have at once:

\proclaim{Proposition 4.8}  Let the hypothesis be as in (4.7).
and let $W$ have a filtration as in (4.7).  Then:
{\pa
\part a $W$ has finite
length as an $\rf$-module if and only if each of the factors
$N_j/L_j$ has finite length  in the category of $R$-modules.
\part b $W$ is \an\ if and only if in some (equivalently,
every filtration) as in (4.6),  the nilpotent factors are all zero.
\par}
\endproclaim

\demb{Proof} (a) This comes down to the assertion that if a power of
$F$ kills $W$ then $W$ has finite length in the category of
$\rf$-modules iff it has finite length in the category
of $R$-modules.  But $W$ has a finite filtration with
factors $\widetilde{F}^j(W)/\widetilde{F}^{j+1}(W)$ on which $F$ acts trivially,
and the result is obvious when $F$ acts trivially. 

It remains to prove (b).  If $W$ is \an, then the
nilpotent factors in any finite filtration must be 0,  since
they are subquotients of $W$.  Now suppose that $W$ has 
a finite filtration by simple $\rf$-modules on which $F$
acts injectively.  Suppose that we have $0 \inc V \inc V' \inc W$
such that $F$ acts trivially on $V'/V$.  This filtration and
the filtration by simple $\rf$-modules on which $F$ acts
injectively have a common refinement in the category
of $\rf$-modules.  This implies that
$V'/V$ has a finite filtration in which  all the factors
are simple $\rf$-modules on which $F$ acts injectively.  Since
$F$ must be zero on the smallest nonzero submodule in the
filtration, this is a contradiction.  
\qed\enddemo

\proclaim{Corollary 4.9} Let $R$ be a local ring of positive
prime characteristic $p$.  
If $M$ is an Artinian $R$-module that is \an\ as
an $\rf$-module, then so is every submodule, quotient module,
and every subquotient of $M$ in the category of 
$\rf$-modules.  \endproclaim
\demb{Proof}  It suffices to show this for submodules and
quotients.  But if $N$ is any $\rf$-submodule, the filtration
$0 \inc N \inc M$  has a common refinement with the filtration
of $M$ with factors that are simple $\rf$-modules with
non-trivial F-action. \qed\enddemo

We also note the following, which is part of [Ly], Theorem (4.2).

\proclaim{Theorem 4.10 (Lyubeznik)} Let $T \to R$ be a surjective
map from a complete regular local ring $T$ of prime characteristic $p >0$
onto a local ring $\rmk$.  Then there exists a contravariant additive
functor $\cH_{T,R}$ from the category of $\rf$-modules that are
Artinian over $R$ to the category of $F_T$-finite modules in the
sense of Lyubeznik such that
{\pa
\part a  $\cH_{T,R}$ is exact. 
\part b \ $\cH_{T,R}(M) = 0$ if and only if the action of some power of
$F$ on $M$ is zero. \par}
\endproclaim

\proclaim{Theorem 4.11} Let $R$ be a local ring of positive
prime characteristic $p$.  Let $M$ be an Artinian $R$-module
that is \an\ as an $\rf$-module.  Then $M$ has only finitely
many $\rf$-submodules. \endproclaim
\demb{Proof} We may replace $R$ by its completion and write $R$ as
$T/J$ where $T$ is a complete regular local ring of characterisitc $p$.
By (4.10) above, we have a contravariant exact functor  $\cH_{T,R}$ 
on $\rf$-modules Artinian over $R$ to $F_T$-finite modules in the
sense of Lyubeznik.  This functor is faithfully exact when restricted
to anti-nilpotent modules, and all subquotients of $M$ are anti-nilpotent.
It follows that if $M_1$ and $M_2$ are distinct $\rf$-submodules of $M$,
then $N_1$ and $N_2$ are distinct, where 
$N_i = \Ker \bigl(\cH(M) \to \cH(M_i)\bigr)$
for $i = 1, 2$.  By the main result of [Ho4],  an $F_T$-finite module in the
sense of Lyubeznik over a regular local ring $T$ has only finitely many 
F-submodules, from which the desired result now follows at once. \qedd

\demb{Discussion 4.12:  local cohomology after adjunction of a formal indeterminate}
Let $R$ be any ring and $M$ an $R$-module.  Let $x$ be a formal power series
indeterminate over  $R$.   We shall denote by $M\inx$ the $R[[x]]$-module
$M \otimes_\Z (\Z[x, \, x^{-1}]/\Z[x])$.  This is evidently an $R[x]$-module, and
since every element is killed by a power of $x$, it is also a module over $R[[x]]$.
Note that if $R$ is an $A$-algebra, this module may also be described as
$M \otimes_A (A[x,\, x^{-1}]/A[x])$.  In particular, $M\inx \cong M \otimes_R (R[x,\,x^{-1}]/R[x])$,
and if $R$ contains a field $K$,  $M\inx \cong M \otimes_K K[x,\,x^{-1}]/K[x]$.  
We have that $M \cong M \otimes x^{-n}$ for all $n \geq 1$ via the map
$u \mapsto  u \otimes x^{-n}$, and we write $Mx^{-n}$ for $M\otimes x^{-n}$. 
As an $R$-module,   $M\inx \cong \bigoplus_{n=1}^\infty Mx^{-n}$,  a countable direct
sum of copies of $M$.  The action of $x$ kills $Mx^{-1}$ and for $n > 1$  maps 
$Mx^{-n} \cong  Mx^{-(n-1)}$
in the obvious way, sending $ux^{-n} \mapsto ux^{-(n-1)}$.  
Then $M \to M\inx$ is a faithfully exact functor from $R$-modules to $R[[x]]$-modules.
If $R$ has prime characteristic $p > 0$ and $M$ is an $\rf$-module,  then we may also
extend this to an $R[[x]] \{F\}$-module structure on $M\inx$  by letting 
$F$ send $ux^{-n} \mapsto F(u)x^{-pn}$.  This gives a convenient way of describing  
what happens to local cohomology when we adjoin a formal power series 
indeterminate to a local ring $R$. \enddemo

\proclaim{Proposition 4.13} Let $R$  be a Noetherian ring, let $I$ be a finitely generated ideal of $R$, 
and let $x$ be a formal power series
indeterminate over $R$. Let $J$ denote the ideal $(I,\,x)R[[x]]$ of $R[[x]]$.
{\pa
\part a For every $i$,  $H^i_J(R[[x]]) \cong H^i_I(R)\inx$. 
 If $R$ has prime characteristic $p >0$
then the action of Frobenius on $H^i_J(R[[x]])$ agrees with the action on $H^i_I(R)\inx$
described above.
\part b  In particular if $\rmk$ is local and and $J = \n$, the maximal ideal of $R[[x]]$,
then for every $i$,  $H^i_\n(R[[x]]) \cong H^i_m(R)\inx$. 
\part c If $\rmk$ is local and $M$ is Artinian, then $M\inx$ is Artinian over $R[[x]]$.  
\part d If $\rmk$ is local of prime characteristic $p > 0$ and $M$ is a simple
$\rf$-module on which the action of $F$ is not 0, then $M\inx$ is a simple $R[[x]]\{F\}$-module.  
\part e If $\rmk$ is local of prime characteristic $p > 0$ and $M$ is an \an\
$\rf$-module, then $M\inx$ is an \an\ $R[[x]]\{F\}$-module.  
\par}
\endproclaim
\bpf (a) Let $\vect f n \in R$ generate $I$. Then  $H^i_I(R)$ is the $i\,$th cohomology
of the complex $C\ub(\uf;\, R)$,  and $H^i_J(R[[x]])$ is the $i\,$th cohomology module 
of the complex $C\ub(\uf,x;\, R[[x]]) \cong C\ub(\uf;\,R) \otimes_R C\ub(x;\, R[[x]])$.
The complex  $C\ub(x;\, R[[x]])$ is simply $$0 \to R[[x]] \to R[[x]]_x \to 0$$ and has
augmentation $R\inx$.  Since $R[[x]]$, $R[[x]]_x$, and $R\inx$ are all
$R$-flat, we have that  $H^i_J(R[[x]])$ is the $i\,$th cohomology module
of the mapping cone of the injection of complexes $C\ub(\uf;\,R[[x]]) \inj C\ub(\uf;\,R[[x]]_x)$,
which may be identified with the cohomology of the quotient complex, and so with
the cohomology of $C\ub(\uf;\,R\inx) \cong C\ub(\uf;\,R) \inx$.  Since $R\inx$ is
$R$-flat (in fact, $R$-free) applying $\blank\otimes_R R\inx$ commutes with formation of 
cohomology, from which the result follows. Part (b) is immediate from part (a).

To prove (c),  note that every element of $M\inx$ is killed by a power of $m$ and of $x$,
and so by a power of $\n$.  It therefore suffices to see that the annihilator of $\n$ is
a finite-dimensional vector space over $K$.  But the annihilator of $x$ is $Mx^{-1}$,
and the annihilator of $m$ in $Mx^{-1}$ is isomorphic with the annihilator of $m$
in $M$.

We next prove (d).  Since the kernel of the action of $F$ on $M$ is an F-stable $R$-submodule
of $M$,  the fact that $M$ is a simple $\rf$-module implies that $F$ acts injectively on $M$.  
Suppose that $N$ is a nonzero $R[[x]]\{F\}$-submodule of $M$, and that 
$u_1x^{-1} + \cdots + u_k x^{-k} \in N$.  By multiplying by $x^{k-1}$ we see that
$u_kx^{-1} \in N$.  Hence, $N$ has nonzero intersection $N_1x^{-1}$ with
$Mx^{-1}$.  $N_1$ is an $R$-submodule of $M$.  It is also F-stable, since if $ux^{-1} \in N$
then $x^{p-1}F(ux^{-1}) = F(u)x^{-1} \in N \cap Mx^{-1}$.  Thus,  $N$ contains
$Mx^{-1}$.     In every degree $h$, let $N_hx^{-h} = N \cap Mx^{-h}$.  Then $N_h \not=0$,
for if $u \in M - \{0\}$ and $q = p^e > h$,  $ x^{q-h}F^e(ux^{-1}) = F^e(u)x^{-h}$,  and $F(u) \not= 0$.
Moreover, the $R$-submodule $N_h \inc M$ is F-stable, for if  $vx^{-h} \in N_hx^{-h}$,
then $x^{ph-h} F(vx^{-h}) = F(v)x^{-h} \in N_hx^{-h}$.  Thus, $N_h = M$ for all $h \geq 1$,
and so $N = M$. 

For part (e),  if $M$ has a finite filtration by simple $\rf$-modules $M_j$ on which $F$ has nonzero
action, then applying $N \mapsto N\inx$  gives a finite filtration of $M\inx$ with factors
$M_j\inx$,  each of which  is a simple $R[[x]]\{F\}$-module by part (d) on which $F$ has
nonzero action.  \qedd

\proclaim{Theorem 4.14} Let $\rmk$ be a local ring of prime charjacteristic $p > 0$.  Let
$x$ be a formal power series indeterminate over $R$.
Let $M$ be an $\rf$-module that is Artinian as an $R$-module.  Then the following are 
equivalent:
{\pa
\part 1  M is \an.
\part 2  $M\inx$  has finite length over $R[[x]]\{F\}$.
\part 3 $M\inx$ has only finitely many F-stable submodules over $R[[x]]$. \par}
When these equivalent conditions hold, $M$ has only finitely many F-stable modules
over $R$. 
\endproclaim
\bpf We show that $(2) \imp (1) \imp (3)$.  Assume (2).   If $M$ is not \an, it has a subquotient
$N \not= 0$ on which the action of $F$ is 0.  Then $N\inx$ is a subquotient of $M\inx$,
and so has finite length as an $R[[x]]\{F\}$-module.  Since $F$ kills it, it must have
finite length as an $R[[x]]$-module.  But this is clearly false, since no power of
$x$ kills $N\inx$.  

To see that $(1) \imp (3)$,  note that by Proposition (4.13e), the fact that   $M$ is \an\ implies
that $M\inx$ is \an\ over $R[[x]]$, and the result now follows from Theorem (4.11).
The implication $(3) \imp (2)$ is obvious. \qedd

The next result is an immediate corollary.

\proclaim{Theorem 4.15} Let $(R,\,m,\,K)$ be a local ring of prime characteristic
$p > 0$ and let $x = x_1$ and $x_2, \, \ldots, \, x_n$ be formal power series 
indeterminates over $R$.  Let $R_n = R[[\vect x n]]$, where $R_0 = R$, and
let $m_n$ be its maximal ideal.   
Then the following conditions on $R$ are equivalent:
{\pa
\part 1  The local cohomology modules $H^i_m(R)$ are anti-nilpotent.
\part 2 The ring  $R[[x]]$ has FH-finite length.
\part 3 The ring $R_n$ is FH-finite for every $n \in \N$.  
\part 4 The local cohomology modules $H^i_{m_n}(R_n)$ are anti-nilpotent over $R_n$ for all $n \in \N$.  
\part 5 The ring $R_n$ has FH-finite length for all $n \in \N$. \par}
When these conditions hold,  $R$ is F-injective.  \endproclaim
\bpf We have that $(1) \imp (4)$ by (4.13a and 4.13e) and a straightforward induction on $n$.
This implies that $R_n$ is FH-finite  for all $n$ by Theorem (4.11).   
Thus, $(4) \imp (3) \imp (5) \imp (2)$, and it suffices
to prove that $(2) \imp (1)$,  which is a consequence of Theorem (4.14).

The statement that $R$ is then F-injective is obvious, since $F$ acts injectively
on any \an\ module. \qedd 

\proclaim{Corollary 4.16} Let $\rmk$ be an F-pure Gorenstein local ring of prime
characteristic $p >0$ and Krull
dimension $d$.   Then $H^d_m(R)$ is \an,  and so $F$ acts injectively on
every subquotient of $H^d_m(R)$.  \endproclaim
\bpf  The hypothesis also holds for $R[[x]]$, and so the result follows from
(3.7) and (4.15) \qedd

\proclaim{Corollary 4.17} Let $\rmk$ be an F-pure Gorenstein local ring prime
characteristic $p >0$ of Krull dimension $d$. Let $J$ be an ideal of $R$ such
that $\di(R/J) = d$.  Then  $H^d_m(R/J)$ is \an, and so $F$ acts injectively
on $H^d_m(R/J)$ (here, $F$ is induced naturally from the Frobenius action $F$ on $R/J$).  Hence, if $R/J$ is Cohen-Macaulay, it is F-injective. \endproclaim

\demb{Proof}  Since $R$ and $R/J$ have the same dimension,
the long exact sequence for local cohomology gives an $\rf$-module surjection
$H^d_m(R) \surj H^d_m(R/J)$,  which shows that $H^d_m(R/J)$ is \an\
as an $\rf$-module, and therefore as an $(R/J)\f$ module as well.
\qed\enddemo  

We shall next need the following result from [Wat]:  the F-pure case is
attributed there to Srinivas.  See Theorem (2.7) of [Wat] and the comment that precedes it.

\proclaim{Theorem 4.18 (Watanabe and Srinivas)}  
Let $h: (R,\,m,\,K) \to (S,\,\n, \, L)$ be a local
homomorphism of local normal domains of prime characterisitic $p > 0$ such
that $S$ is module-finite over $R$ and the map $h$ is \'etale in codimension one.
If $R$ is strongly F-regular, then so is $S$.  If $R$ is F-pure, then so is $S$.  
\endproclaim

The explicit statement in [Wat] is for the F-regular case, by  which the author
means the weakly F-regular case.  However, the proof given uses the
criterion (i) of [Wat], Proposition (1.4), that the local ring $\rmk$ is weakly F-regular
if and only if $0$ is tightly closed in $E_R(K)$, 
which is correct for finitistic tight closure
but not for the version of tight closure being used in [Wat].  In fact, condition
(i) as used in [Wat] characterizes strong F-regularity in the F-finite case, and 
we take it as the definition of  strong F-regularity here.   

On the other hand, 
there are no problems whatsoever in proving the final statement about 
F-purity.  The action of Frobenius
$F_S:E_S(L) \to F_S\bigl(E_S(L)\bigr)$ is shown to be the same as the
action of Frobenius when $E_S(L)$ is viewed as $R$-module.  Since
$R$ is F-pure, the Frobenius action $M \to F_R(M)$ is injective for
any $R$-module.

\proclaim{Corollary 4.19 (Watanabe)} Let $(R,\,m)$ be a normal local domain of 
characteristic $p > 0$. Let $I$ be an ideal of pure height one, and suppose that
$I$ has finite order $k > 1$ in the divisor class group of $R$.  Choose a generator
$u$ of  $I^{(k)}$.   We let 
$$S = R \oplus It \oplus \cdots \oplus I^{(j)}\oplus \cdots \oplus I^{(k-1)}$$
with $I^{(k)}$ identified with $R$ using the isomorphism $R \cong I^{(k)}$ such that 
$1 \mapsto u$.   (If  $t$ is an indeterminate, we can give the following more formal 
description:  form $T = \bigoplus_{j=0}^\infty I^{(j)} t^j \inc R[t]$, and let  $S = T/(ut^k - 1$.)
This ring is module-finite over $R$, and if $k$ is relatively prime to $p$, it is 
\'etale over $R$ in codimension one.  

Hence, if  $k$ is relatively prime to $p$, then $S$ is strongly F-regular
if and only if  $R$ is strongly F-regular, and $S$ is F-pure if and only if 
$R$ is F-pure.

Moreover, if  $I \cong \omega$ is
a canonical module for $R$,  then $S$ is quasi-Gorenstein.
\endproclaim 

The final statement is expected because, by the discussion of 
canonical modules for module-finite extensions in (2.2),  we have
that $\omega_S \cong \Hom_R(S, \, \omega)$ and $\Hom_R(\omega^{(i)},\, \omega)
\cong \omega^{-(i-1)} \cong \omega^{(k-(i-1))}$. 
See [Wa], [TW], and [YW, \S 3], and [Si] for further details and background 
on this technique.  The result above will enable us to use 
this ``canonical cover trick" to prove the Theorem below by 
reduction to the quasi-Gorenstein case.  A word of caution is in order:
even if $R$ is Cohen-Macaulay, examples in [Si] show that the
auxiliary ring $S$ described in (4.19) need not be Cohen-Macaulay,
and so one is forced to consider the quasi-Gorenstein property.  There
are examples (see [Si, Theorem (6.1)]) where $R$ is F-rational but 
$S$ is not Cohen-Macaulay.
On the other hand, if $R$ is strongly F-regular, the result of [Wat] shows
that $S$ is as well:  in particular, $S$ is Cohen-Macaulay in this case.

\proclaim{Theorem 4.20} Let $R$ be a Cohen-Macaulay F-pure normal local domain
of Krull dimension $d$
such that $R$ has canonical module $\omega = \omega_R$ of finite order
$k$ relatively prime to $p$ in the divisor class group of $R$.  Then $H_m^d(R)$
is anti-nilpotent, so that $R$ is FH-finite.
\endproclaim
\bpf Since the hypotheses are stable under adjunction of a power series
indeterminate, it follows from Theorem (4.15) that it is sufficient to show that
$R$ is $FH$-finite.  We may identify $\omega$ with a pure height one ideal
of $I$ of $R$.  We form the ring $S$ described in (4.19).  Then $S$ is
F-pure and quasi-Gorenstein, and so $H^d_m(S)$ has only finitely many
F-stable submodules by  Theorem (3.7).  The same holds for $H^d_m(R)$ by
Lemma (2.7b), while the lower local cohomology modules of $R$ with
support in $m$ vanish. \qedd

The following improves the conclusion of (2.14) with some additional
hypotheses. 

\proclaim{Theorem 4.21} Let $\rmk$ be an F-injective Cohen-Macaulay 
local ring with of prime characteristic $p>0$ such that $K$ is
perfect.  Suppose that $R$ has an $m$-primary ideal $\fA$ such that
$\fA I^* \inc I$ for every ideal $I$ generated by a system of parameters.
Let $d = \di(R)$.  Then $H^d_m(R)$ is anti-nilpotent,
so that $R$ and every formal power series ring over $R$ is FH-finite.
\endproclaim
\bpf  Let $H = H^d_m(R)$ and let $V = 0^*_H$,  which, as in the proof
of (2.14), is killed by $\fA$ and has finite length.  Then $R/\fA$ is
a complete local ring with a perfect residue class field, and contains
a unique coefficient field $K$.  This gives $V$ the structure of a $K$-module,
i.e., it is a finite-dimensional $K$-vector space,  and $F:V \to V$ is $K$-linear
if we let the action of $K$ on the second copy of $V$ be such that $c \cdot v = c^pv$
for $c \in K$.  Since $K$ is perfect, the dimension of $V$ does not change when
we restrict $F$ in this way.  Since $R$ is F-injective, the action of $F$ on $V$ is
then a vector space isomorphism, and is then also an isomorphism when restricted
to subquotients that are $K\{F\}$-modules.  It follows that $V$ is anti-nilpotent
over $\rf$,  and to complete the proof it will suffice to show that $F$ cannot
act trivially on the the simple $\rf$-module $H/V$.  

Choose a system of 
parameters $\vect x d$ for $R$. Let $I_t = (x_1^t, \, \ldots, \, x_d^t)R$.
For any sufficiently large value of $t$, we may identify $V$ with 
$I^*/I$.  If  $F$ acts trivially on $H/V$,  then for all large $t$, the image of
$1 \in R/I_t \inc H$ under $F$ is 0 in $H/V$,  which means that $(x_1 \cdots x_d)^{p-1}
\in I_{pt}^*$,  and then $\fA (x_1 \cdots x_d)^{p-1} \inc I_{pt}$.  This implies
that  $$\fA \inc I_{pt}:_R(x_1 \cdots x_d)^{p-1} = I_{pt-p+1}$$
for all $t \gg 0$, which is clearly false.  \qedd

\quad\bigskip\bigskip
{\chbf 
\centerline{\hdbf 5. FACE RINGS}
}
\bigskip
We give a brief treatment of the decomposition of the local cohomology
of face rings over a field with support in the homogeneous maximal ideal.
This is discussed in [BH], \S5.3, although not in quite sharp enough a
form for our needs here, and there are sharp results in substantially
greater generality in [BBR]: see Theorem (5.5), p.\ 218.  However, neither
result discusses the $\rf$-structure.  

Let $K$ be a fixed field of positive characteristic $p$ and
let $\Delta$ be an abstract finite simplicial complex with
vertices $\vect x n$. Let $I_\dl$ denote the ideal in the
polynomial ring $S = K[\vect x n]$ generated by all monomials
in the $x_j$ 
such that the set of variables occurring in the monomial is
not a face of $\dl$.  This ideal is evidently generated by the square-free
monomials in the $x_j$ corresponding to minimal subsets of the
variables that are not faces of $\dl$.  Let $K[\dl] = S/I_\dl$,
the {\it face ring} (or {\it Stanley-Reisner ring}) of $\dl$ over $K$.
The minimal primes of $K[\dl]$ correspond to the maximal faces $\sigma$
of $\dl$:  the quotient by the minimal prime corresponding to 
$\sigma$ is a polynomial ring in the variables occurring in $\sigma$.
The Krull dimension of $K[\dl]$ is therefore one more than
the dimension of the simplicial complex $\dl$.   

If $\sigma$ is any face of $\dl$,  the link, denoted $\link(\sigma)$, of $\sigma$
in $\dl$ is the abstract simplicial complex consisting of all
faces $\tau$ of $\dl$ disjoint from $\sigma$ such that
$\sigma \cup \tau \in \dl$.   The link of the empty face is 
$\dl$ itself.   By a theorem of G.\ Reisner [Rei],  
$K[\dl]$ is Cohen-Macaulay
if and only if the reduced simplicial cohomology of every link
vanishes except possibly in the top dimension, i.e., in the dimension
of the link itself.  

Note that the reduced simplicial cohomology $\wt{H}^i(\dl;\, K)$ of a finite simplicial
complex $\dl \not= \emptyset$ is the same as the simplicial comhomology 
unless $i = 0$,  in which
case its dimension as a $K$-vector space is one smaller.  If $\dl$ is an
$i$-simplex, the reduced simplicial cohomology  vanishes in all dimensions,
unless $\dl$ is empty,  in which case we have $\wt{H}^{-1}(\emptyset;\,K) \cong K$:
$\wt{H}^i(\emptyset;\,K) = 0$ for all $i \not= -1$.
Note also that $\emptyset$ is the only simplicial complex $\dl$ such that
$\wt{H}^i(\dl;\,K) \not=0$ for a value of $i < 0$.   

Let $m$ be the homogeneous maximal ideal of $K[\dl]$.  We shall show
that $K[\dl]_m$ and its completion are FH-finite in all cases, and in fact,
the local cohomology modules are anti-nilpotent.  This follows from
the following theorem, which also recovers Reisner's result [Rei] mentioned above in 
a finer form:  it also gives a completely explicit description of all the $H^i_m(K[\dl])$,
including their structure as $\rf$-modules.  We  write $|\nu|$ for the cardinality
of the set $\nu$.   If $\nu \in \dl$,  we let $$K[\nu] = K[\dl]/(x_i: x_i \notin \nu),$$
which is a $K[\dl]$-algebra and is also the polynomial ring over $K$ in the
variables $x_j$ that are vertices of $\nu$.  Then $H^i_m(S_\nu)$ vanishes
except when $i  = |\nu|$. When $i = \nu$ it is the highest nonvanishing
local cohomology of a polynomial ring, and, if the characteristic of $K$ is
$p > 0$, it is a simple $\rf$-module on which  $F$ acts injectively.

Note that if $p>0$ is prime, $\kappa = \Z/p\Z$, 
$R$ and $K$ are rings of characteristic $p$, and $H$ is an 
$\rf$-module, $K \otimes_{\kappa}H$
has the structure of a $(K \otimes_\kappa R)\f$-module:
the action of $F$ is determined by the rule 
$F(c \otimes u) = c^p \otimes F(u)$ for all $c \in K$ and $u \in H$.  
This is well-defined because the action of $F$ restricts to the 
identity map on $\Z/p\Z$.  

\proclaim{Theorem 5.1} With $R = K[\dl]$ as above, let $\kappa$ denote
the prime field of $K$. Let $m$ and $\mu$ be the homogeneous  maximal ideals 
of $R$ and $\kappa[\dl]$, respectively.  Then
$$ (*) \quad H^i_m(R) \cong  \bigoplus_{\nu \in \dl} \wt{H}^{i-1 - |\nu|}(\link(\nu);\, K) 
\otimes_\kappa H^{|\nu|}_{\mu}(\kappa[\nu]).$$
If $K$ has characterisitc $p > 0$, this is also an isomorphism of $\rf$-modules,
with the action of $F$ described in the paragraph above.
Hence, every $H^i_m(R)$ is a finite direct sum of simple $\rf$-modules on which
$F$ acts injectively.  

If $(R_1,\,m_1)$ is either $R_m$ or its completion,  then for all $i$, 
$H^i_m(R)$  may be identified with $H^i_{m_1}(R_1)$,
and $H^i_{m_1}(R_1)$ is a finite
direct sum of simple $R_1\f$-modules on which $F$ acts injectively, and
so is anti-nilpotent and FH-finite over $R_1$. \endproclaim

\bpf   If  $\sigma$ is a subset of the 
$x_j$ we denote by $x(\sigma)$
the product of the $x_j$ for $j \in \sigma$.  Thus, in $K[\dl]$,
the image of $x(\sigma)$ is nonzero if and only if $\sigma \in \dl$. 
Our convention is that $\sigma = \emptyset$ is in $\dl$,  is the
unique face of dimension $-1$, and that $x(\emptyset) = 1$.  
We write $[\dl]_i$ for the set of faces of $\dl$ of dimension $i$. 
Then $H^i_m(K[\dl])$ is the $i\,$th cohomology module of the
complex $C\ub$ whose $i\,$th term is displayed below:
$$
0 \to K[\dl] \to \cdots \to \bigoplus_{\sigma \in \dl_{i-1}}K[\dl]_{x(\sigma)} \to \cdots 
$$
The initial nonzero term $K[\dl]$ may be thought of as $K[\dl]_{x(\emptyset)}$
and the highest nonzero terms occur in degree $\dim(\dl)+1$.  
Let $\theta = (\vect \theta n) \in \Z^n$ by an $n$-tuple of integers.
We want to calculate that $\theta$-graded piece of $H^i_m(K[\dl])$.  
This is the same as the $i\,$th cohomology of the $\theta$-graded
piece of the complex:  denote by $C^\bullet[\theta]$ the $\theta$-graded piece of
the complex $C\ub$.  Let $\ng(\theta)$ (respectively, $\pos(\theta)$,
respectively, $\spu(\theta)\,$) 
denote the set of variables $x_i$ such that $\theta_i$ is stictly negative (respectively,
strictly positive, respectively, nonzero). Thus, $\spu(\theta)$ is the disjoint
union of $\ng(\theta)$ and $\pos(\theta)$.     

Let $\nu = \ng(\theta)$ and $\pi = \pos(\theta)$.  
Then $K[\dl]_{x(\sigma)}$ has a nonzero 
component in degree $\theta$ if and only if  $\nu \inc \sigma$ and 
$\sigma \cup \pi \in \dl$,
and then there is a unique copy of $K$ corresponding to $\theta$
in the complex.   By deleting the variables occurring in
$\nu = \ng(\theta)$ from each face,  we find that $C^\bullet[\theta]$ 
corresponds, with a shift in degree by the cardinality $|\nu|$ of $\nu$, 
to the complex used to calculate the reduced
simplicial cohomology of $\dl_{\nu,\pi}$,  where this is the
subcomplex of the link of $\nu$ consisting of all simplices  $\tau$
such that $\tau \cup \pi \in \dl$.  If $\pi$ is non-empty,  
$\dl_{\nu, \pi}$ is a cone on any vertex in $\pi$.  Hence, the graded
component of local cohomology in degree $\theta$ can be nonzero
only when $\pos(\theta) =\emptyset$ and $\nu = \ng(\theta) \in \dl$.
It now follows that
$$
[H^i_m(K[\dl])]_\theta \cong \bigoplus_{\nu} \wt{H}^{i - |\nu| -1}(\hbox{link}\,\nu)x^{\theta}
$$
if $\pi = \emptyset$ and $\nu \in \dl$ is the set of variables corresponding
to strictly negative entries in $\theta$, and is zero otherwise. 

It follows that we may identify
$$
H^i_m(K[\dl]) \cong 
\bigoplus_{\nu \in \dl}\quad 
\bigl(\bigoplus_{\spus (w) = \ngs (w) = \nu} \wt{H}^{i-|\nu|-1}(\hbox{link}\,\nu)w\bigr) ,  
$$
where $w$ runs through all monomials with non-positive exponents such that the
set of variables with strictly negative exponents is $\nu$. 

We next want to show that if $\nu \in \dl$,  then
$$\bigoplus_{\spus (\theta) = \ngs(\theta) = \nu} [H^i_m(K[\dl]]_\theta 
\cong \wt{H}^{i-|\nu|-1}(\link(\nu);\, K) \otimes_K H^{|\nu|}_m(K[\nu]).$$
The term on the right may also be written as 
$\wt{H}^{i-|\nu|-1}(\link(\nu);\, K) \otimes_\kappa H^{|\nu|}_\mu(\kappa[\nu])$.
We also need to check that the actions of $F$ agree.

The action of $F$ on $C\ub(\ux;\, K[\dl])$ is obtained from the action of $F$
on $C\ub(\ux;\,\kappa[\dl]) $ by applying $K \otimes_{\kappa} \blank$. Thus, we reduce
at once to the case where $K = \kappa$,  which we assume henceforth.  

Let  $\gamma \otimes w$  be an element in the cohomology, where
$\gamma \in  \wt{H}^{i-|\nu|-1}(\link(\nu);\, K)$ and  $w$ is a monomial 
in the variables of $\nu$ with
all exponents strictly negative.   The action of $x_i$ by multiplication is obvious
in most cases.  If  $x_i \notin  \nu$, the product is 0.  If  $x_i$ occurs with
an exponent other than $-1$ in $w$,  one simply gets  $\gamma \otimes (x_iw)$. 
The main non-trivial  point is that  if $x_i$ occurs with exponent $-1$ in $w$,
$x_i$ kills $\gamma \otimes w$.  To verify this, let $\nu' = \nu-\{x_i\}$. 
Take a cocycle $\eta$ that
represents  $\gamma$.    After we multiply
by $x_i$,  we get an element of $H^{i-|\nu|}\bigl(\link(\nu')\bigr)$. Note that each
simplex remaining when we delete the variables in $\nu'$ involves $x_i$,
and so that the cocycle $\eta'$ we get from $\eta$ may be viewed as a cocycle 
of the complex used to compute the reduced simplicial cohomology of the
closed star of $x_i$ in $\link(\nu')$. Since this closed star is a cone, that cohomology
is 0.  This shows that $x_i$ kills every homogeneous component whose 
degree has $-1$ in the $i\,$th coordinate.

We have completed the calculation of the structure of the local cohomology as
an $R$-module.  On the other hand, given $\nu \in \dl$, 
because the field is $\kappa$,  when
$F$ acts on the complex 
$$\bigoplus_{\spus(\theta) = \neg(\theta) = \nu} [C\ub(\ux;\,R)]_\theta$$
the value of $F$ acting on an element of the form  $\eta w$,  where $\eta$ is a 
cocycle, if simply $\eta w^p$, and so it follows that 
$$F(\gamma \otimes w) = \gamma \otimes w^p.$$  This shows that $\rf$-module 
structure is preserved by the isomorphism $(*)$.  \qedd

} %% End baselineskip  %%%  REINSERT WHEN FINISHED
%%\vfill\eject %% may want to change this
\quad

{
\smc
\baselineskip = 10 pt
\settabs 7 \columns
\quad\bigskip
\+ Department of Mathematics                     &&&& Department of Mathematics\cr
\+ \ \ \ \ \ \ \ \ \ \ $\,$and Statistics                            &&&& University of Michigan\cr
\+ Georgia State University                          &&&& East Hall, 530 Church St.\cr
\+  750 COE, 7th floor,\ 30 Pryor St.            &&&& Ann Arbor, MI 48109--1043\cr
\+  Atlanta, GA 30303--3083                        &&&& USA\cr
\+ USA   \cr                      
\smallskip
\+E-mail:                        &&&&E-mail:\cr
\vskip 1.5 pt plus .5 pt minus .5 pt
{\rm
\+fenescu\@gsu.edu                              &&&&hochster\@umich.edu\cr
}
}

\enddocument